\documentstyle[11pt]{article}
  \newcommand{\bref}[5]{\noindent\parbox[t]{0.8cm}{#1}\parbox[t]{15.2cm}{{#2}{\it #3}{\bf #4}#5}\par}
  
  \def\sg{\sigma}

\begin{document}
  \title{\bf {Exact solution to an extremal problem on
graphic sequences with a realization containing every $2$-tree on
$k$ vertices}
  \thanks{Supported by Hainan Provincial Natural Science Foundation of China
(No. 118QN252) and National Natural Science Foundation of China (No.
11561017).}}
\author{De-Yan Zeng$^{1,2}$, Dong-Yang Zhai$^2$, Jian-Hua Yin$^1$\thanks{Corresponding author.\ \ E-mail: yinjh@hainu.edu.cn}\\
{\small $^1$Department of Mathematics, College of Information
Science and Technology,}\\
{\small Hainan University, Haikou 570228, P.R. China}\\
{\small $^2$Institute of Science and Technology, University of
Sanya, Sanya 572022, P.R. China.}}
\date{}
\maketitle

\begin{center}
\parbox{0.9\hsize}
{\small {\bf Abstract.}\ \ A simple graph $G$ is an {\it 2-tree} if
$G=K_3$, or $G$ has a vertex $v$ of degree 2, whose neighbors are
adjacent, and $G-v$ is an 2-tree. Clearly, if $G$ is an 2-tree on
$n$ vertices, then $|E(G)|=2n-3$. A non-increasing sequence
$\pi=(d_1,\ldots,d_n)$ of nonnegative integers is a {\it graphic
sequence} if it is realizable by a simple graph $G$ on $n$ vertices.
Yin and Li (Acta Mathematica Sinica, English Series,
25(2009)795--802) proved that if $k\ge 2$, $n\ge
\frac{9}{2}k^2+\frac{19}{2}k$ and $\pi=(d_1,\ldots,d_n)$ is a
graphic sequence with $\sum\limits_{i=1}^n d_i>(k-2)n$, then $\pi$
has a realization containing every 1-tree (the usual tree) on $k$
vertices. Moreover, the lower bound $(k-2)n$ is the best possible.
This is a variation of a conjecture due to Erd\H{o}s and S\'{o}s. In
this paper, we investigate an analogue problem for $2$-trees and
prove that if $k\ge 3$ is an integer with $k\equiv i(\mbox{mod }3)$,
$n\geq20\lfloor\frac{k}{3}\rfloor^2+31\lfloor\frac{k}{3}\rfloor+12$
and $\pi=(d_1,\ldots,d_n)$ is a graphic sequence with
$\sum\limits_{i=1}^n
d_i>\max\{(k-1)(n-1),2\lfloor\frac{2k}{3}\rfloor
n-2n-\lfloor\frac{2k}{3}\rfloor^2+\lfloor\frac{2k}{3}\rfloor+1-(-1)^i\}$,
then $\pi$ has a realization containing every 2-tree on $k$
vertices. Moreover, the lower bound
$\max\{(k-1)(n-1),2\lfloor\frac{2k}{3}\rfloor
n-2n-\lfloor\frac{2k}{3}\rfloor^2+\lfloor\frac{2k}{3}\rfloor+1-(-1)^i\}$
is the best possible. This result implies a conjecture due to Zeng
and Yin (Discrete Math. Theor. Comput. Sci., 17(3)(2016), 315--326).
}
\end{center}

\section*{1. Introduction}
\hskip\parindent  Let $K_m$, $K_{m,n}$ and $P_k$ denote the complete
graph on $m$ vertices, the $m\times n$ complete bipartite graph and
the path on $k$ vertices, respectively. Terms and notation not
defined here are from [1]. A simple graph $G$ is an {\it 2-tree} if
$G=K_3$, or $G$ has a vertex $v$ of degree 2, whose neighbors are
adjacent, and $G-v$ is an $2$-tree. It is easy to see that if $G$ is
an 2-tree on $n$ vertices, then $|E(G)|=2n-3$. An {\it ear} in an
2-tree is a vertex of degree 2 whose neighbors are adjacent.

The set of all non-increasing sequences $\pi=(d_1,\ldots,d_n)$ of
nonnegative integers with $d_1\le n-1$ is denoted by $NS_n$. A
sequence $\pi\in NS_n$ is said to be {\it graphic} if it is the
degree sequence of a simple graph $G$ on $n$ vertices, and such a
graph $G$ is called a {\it realization} of $\pi$. The set of all
graphic sequences in $NS_n$ is denoted by $GS_n$. For a sequence
$\pi=(d_1,\ldots,d_n)$, we denote $\sigma(\pi)=d_1+\cdots+d_n$.
 Yin and Li [14] investigated a variation of a conjecture due to
Erd\H{o}s and S\'{o}s (see [1], Problem 12 in page 247), that is, an
extremal problem for a sequence $\pi\in GS_n$ to have a realization
containing every 1-tree (the usual tree) on $k$ vertices as a
subgraph. They proved the following Theorem 1.1.

{\bf Theorem 1.1} [14]\ \ {\it If $k\ge 2$, $n\ge
\frac{9}{2}k^2+\frac{19}{2}k$ and $\pi=(d_1,\ldots,d_n)\in GS_n$
with $\sg(\pi)>(k-2)n$, then $\pi$ has a realization containing
every 1-tree on $k$ vertices. Moreover, the lower bound $(k-2)n$ is
the best possible.}

This kind of extremal problem was firstly introduced by Erd\H{o}s et
al. (see [5--6]). Zeng and Yin [15] investigated an analogous
extremal problem for a sequence $\pi\in GS_n$ to have a realization
containing every 2-tree on $k$ vertices as a subgraph. They
established the following Theorem 1.2--1.3.

{\bf Theorem 1.2} [15]\ \ {\it If $k\geq 3$, $n\ge 2k^2-k$ and
 $\pi=(d_1,\ldots,d_n)\in GS_n$ with $\sigma(\pi)>\frac{4kn}{3}
 -\frac{5n}{3}$, then $\pi$ has a realization containing every 2-tree
  on $k$ vertices as a subgraph.}

{\bf Theorem 1.3} [15]\ \ {\it For $k\ge 3$ with $k\equiv
i(\mbox{mod }3)$, there exists a sequence $\pi\in GS_n$ with
$\sigma(\pi)=2\lfloor\frac{2k}{3}\rfloor
n-2n-\lfloor\frac{2k}{3}\rfloor^2+\lfloor\frac{2k}{3}\rfloor+1-(-1)^i$
such that $\pi$ has no realization containing every 2-tree on $k$
vertices.}


For $k\ge 3$ with $k\equiv i(\mbox{mod }3)$, Zeng and Yin [15] felt
that $2\lfloor\frac{2k}{3}\rfloor
n-2n-\lfloor\frac{2k}{3}\rfloor^2+\lfloor\frac{2k}{3}\rfloor+1-(-1)^i$
is the best possible lower bound for sufficiently large $n$, thus
they proposed the following conjecture 1.1.

{\bf Conjecture 1.1} [15]\ \ {\it If $k\ge 3$ with $k\equiv
i(\mbox{mod }3)$, $n$ is sufficiently large, and $\pi\in GS_n$ with
$\sigma(\pi)> 2\lfloor\frac{2k}{3}\rfloor
n-2n-\lfloor\frac{2k}{3}\rfloor^2+\lfloor\frac{2k}{3}\rfloor+1-(-1)^i$,
then $\pi$ has a realization containing every 2-tree on $k$
vertices. Moreover, the lower bound $2\lfloor\frac{2k}{3}\rfloor
n-2n-\lfloor\frac{2k}{3}\rfloor^2+\lfloor\frac{2k}{3}\rfloor+1-(-1)^i$
is the best possible.}

In this paper, we further obtain the following Theorem 1.4--1.5. For
convenience, we denote $N(3)=6,N(4)=7,N(5)=24,N(7)=93$ and for $k=6$
or $k\ge 8$,
$$N(k)=\left\{\begin{array}{ll}
20\lfloor\frac{k}{3}\rfloor^2-\lfloor\frac{k}{3}\rfloor,& \mbox{if $k\equiv 0(mod\ 3)$,}\\
20\lfloor\frac{k}{3}\rfloor^2+23\lfloor\frac{k}{3}\rfloor+5,& \mbox{if $k\equiv 1(mod\ 3)$,}\\
20\lfloor\frac{k}{3}\rfloor^2+31\lfloor\frac{k}{3}\rfloor+12,& \mbox{if $k\equiv 2(mod\ 3)$.}
\end{array}\right.$$

{\bf Theorem 1.4} \ \ {\it If $k\in\{3,4,5,7\}$, $n\geq N(k)$ and
$\pi=(d_1,\ldots,d_n)\in GS_n$ with $\sigma(\pi)>(k-1)(n-1)$, then
$\pi$ has a realization containing every 2-tree on $k$ vertices.
Moreover, the lower bound $(k-1)(n-1)$ is the best possible.}

{\bf Theorem 1.5} \ \ {\it If $k=6$ or $k\ge 8$ with $k\equiv
i(\mbox{mod }3)$, $n\ge N(k)$ and
 $\pi=(d_1,\ldots,d_n)\in GS_n$ with $\sigma(\pi)>2\lfloor\frac{2k}{3}\rfloor
n-2n-\lfloor\frac{2k}{3}\rfloor^2+\lfloor\frac{2k}{3}\rfloor+1-(-1)^i$,
then $\pi$ has a realization containing every 2-tree on $k$
vertices. Moreover, the lower bound $2\lfloor\frac{2k}{3}\rfloor
n-2n-\lfloor\frac{2k}{3}\rfloor^2+\lfloor\frac{2k}{3}\rfloor+1-(-1)^i$
is the best possible.}

Theorem 1.5 implies that Conjecture 1.1 is true for $k=6$ or $k\ge
8$.

\section*{2. Useful Known Results}

\hskip\parindent In order to prove Theorem 1.4--1.5, we need some
known results. Let $\pi=(d_1,\ldots,d_n)\in NS_n$ and $k$ be an
integer with $1\le k\le n$. Let
$$\pi_k''=\left\{\begin{array}{ll}
(d_1-1,\ldots,d_{k-1}-1,d_{k+1}-1,\ldots,d_{d_k+1}-1,d_{d_k+2},\ldots,d_n),&
\mbox{if $d_k\ge k$,}\\
(d_1-1,\ldots,d_{d_k}-1,d_{d_k+1},\ldots,d_{k-1},d_{k+1},\ldots,d_n),
& \mbox{if $d_k<k$.}
\end{array}\right.$$
Let $\pi_k'=(d_1',\ldots,d_{n-1}')$, where $d_1'\ge \cdots\ge
d_{n-1}'$ is a rearrangement in non-increasing order of the $n-1$
terms of $\pi_{k}''$. We say that $\pi_k'$ is the {\it residual
sequence} obtained from $\pi$ by laying off $d_k$. It is easy to see
that if $\pi_k'$ is graphic then so is $\pi$, since a realization
$G$ of $\pi$ can be obtained from a realization $G'$ of $\pi_k'$ by
adding a new vertex of degree $d_k$ and joining it to the vertices
whose degrees are reduced by one in going from $\pi$ to $\pi_k'$. In
fact, more is true:

{\bf Theorem 2.1} [7]\ \ {\it $\pi\in GS_n$ if and only if
$\pi_k'\in GS_{n-1}$.}

{\bf Theorem 2.2} [4]\ \ {\it Let $\pi=(d_1,\ldots,d_n)\in NS_n$,
where $\sigma(\pi)$ is even. Then $\pi \in GS_n$ if and only if
$\sum\limits_{i=1}^{t}d_i\leq
t(t-1)+\sum\limits_{i=t+1}^{n}\min\{t,d_i\}$ for each $t$ with
$1\leq t\leq n-1$.}

{\bf Theorem 2.3} [13]\ \ {\it Let $\pi=(d_1,\ldots,d_n)\in NS_n$,
where $d_1=m$ and $\sigma(\pi)$ is even. If there exist $n_1\leq n$
and $h\geq1$ such that $d_{n_1}\geq h$ and
$n_1\geq\frac{1}{h}\lfloor\frac{(m+h+1)^2}{4}\rfloor$, then $\pi\in
GS_n$.}

{\bf Theorem 2.4} [6]\ \ {\it If $\pi=(d_1,\ldots,d_n)\in NS_n$ has
a realization $G$ containing $H$ as a subgraph, then there exists a
realization $G'$ of $\pi$ containing $H$ on those vertices with
degrees $d_1,\ldots,d_{|V(H)|}$.}

{\bf Theorem 2.5} [12]\ \ {\it Let $n\geq r$ and
$\pi=(d_1,\ldots,d_n)\in GS_n$ with $d_r\geq r-1$. If $d_i\geq
2r-2-i$ for $i=1,\ldots,r-2$, then $\pi$ has a realization
containing $K_r$.}

{\bf Theorem 2.6} [11]\ \ {\it If $r\ge 3$, $n\geq2r-1$ and
$\pi=(d_1,\ldots,d_n)\in GS_n$ with $\sigma(\pi)\ge 2n(r-2)+2$, then
$\pi$ has a realization containing $K_r$.}

{\bf Theorem 2.7}\ \ {\it Let $\pi=(d_1,\ldots,d_n)\in GS_n$.

(1) [5]\ \ If $n\ge 6$ and $\sg(\pi)> 2n-2$, then $\pi$ has a
realization containing $K_3$.

(2) [8]\ \ If $n\ge 7$ and $\sg(\pi)> 3n-3$, then $\pi$ has a
realization containing $K_4-e$, where $K_4-e$ is the graph obtained
from $K_4$ by removing an edge.

(3) [9]\ \ If $n\ge 24$ and $\sg(\pi)> 4n-4$, then $\pi$ has a
realization containing $K_5-E(P_3)$, where $K_5-E(P_3)$ is the graph
obtained from $K_5$ by removing all edges of a path $P_3$ in $K_5$.

(4) [10]\ \ If $n\ge 5$ and $\sg(\pi)> 4n-6$, then $\pi$ has a
realization containing $K_5-E(P_4)$, where $K_5-E(P_4)$ is the graph
obtained from $K_5$ by removing all edges of a path $P_4$ in $K_5$.
}

We note that an 2-tree can be constructed from an edge by repeatedly
adding a new vertex and making it adjacent to the two ends of an
edge in the graph formed so far. We refer to the initial edge in
constructing such an 2-tree as a {\it base} of the 2-tree. Some
properties of 2-trees can be summarized as follows.

{\bf Theorem 2.8} [2,3]\ \ {\it Let $G$ be any 2-tree with $n\geq3$
vertices. Then

(1)\ \ $G$ has at least two ears;

(2)\ \ Every vertex of degree 2 in $G$ is an ear;

(3)\ \ No two ears in $G$ are adjacent unless $G=K_3$;

(4)\ \ $G$ does not contain any chordless cycle of length at least
4;

(5)\ \ $G$ is 2-connected;

(6)\ \ Every edge of $G$ can be a base.}

\section*{3. Proof of Theorem 1.4}

\hskip\parindent  We know that $G$ is an 2-tree if either $G=K_3$,
or $G$ has an ear $u$ such that $G'=G-u$ is an 2-tree. In order
words, every 2-tree $G\not=K_3$ can be obtained from some 2-tree
$G'$ by adding a new vertex $u$ adjacent to two vertices, $v$ and
$w$, where $vw\in E(G')$. We call this process {\it attaching} $u$
to $vw$ and denote $vw=e(u)$. Hence the all 2-trees on 7 vertices
can be gotten
in Figure 1. \\[2mm]

\unitlength=1mm
\begin{picture}(150,25)
\put(7.5,12.5){\circle*{1}}         \put(15,12.5){\circle*{1}}
\put(22.5,12.5){\circle*{1}}         \put(37.5,12.5){\circle*{1}}
\put(45,12.5){\circle*{1}}       \put(30,5){\circle*{1}}
\put(30,20){\circle*{1}}         \put(30,5){\line(-3,1){22.5}}
\put(30,5){\line(-2,1){15}}      \put(30,5){\line(-1,1){7.5}}
\put(30,5){\line(0,2){15}}       \put(30,5){\line(1,1){7.5}}
\put(30,5){\line(2,1){15}}       \put(30,20){\line(-3,-1){22.5}}
\put(30,20){\line(-2,-1){15}}    \put(30,20){\line(-1,-1){7.5}}
\put(30,20){\line(1,-1){7.5}}    \put(30,20){\line(2,-1){15}}
\put(2.5,11.5){$u_3$}\put(17,11.5){$u_4$}\put(24.5,11.5){$u_5$}
\put(38.5,11.5){$u_6$}\put(47,11.5){$u_7$}\put(23,1){$u_1(y,v)$}
\put(23,22){$u_2(x,u)$}\put(28,-4){$T_1$}
\put(73,-4){$T_2$}\put(118,-4){$T_3$}

\put(60,5){\circle*{1}}            \put(58,1){$u_7$}
\put(75,5){\circle*{1}}            \put(68,1){$u_2(y,v)$}
\put(90,5){\circle*{1}}            \put(88,1){$u_5$}
\put(60,20){\circle*{1}}            \put(58,22){$u_4$}
\put(75,20){\circle*{1}}           \put(68,22){$u_3(x,u)$}
\put(90,20){\circle*{1}}           \put(88,22){$u_6$}
\put(82.5,12.5){\circle*{1}}      \put(84.5,11.5){$u_1(w)$}
\put(105,12.5){\circle*{1}}         \put(100,11.5){$u_7$}
\put(112.5,12.5){\circle*{1}}       \put(107.5,11.5){$u_4$}
\put(120,5){\circle*{1}}           \put(111,1){$u_1(y,v,w)$}
\put(135,5){\circle*{1}}           \put(133,1){$u_5$}
\put(120,20){\circle*{1}}          \put(111,22){$u_3(x,u)$}
\put(127.5,20){\circle*{1}}        \put(125.5,22){$u_6$}
\put(135,20){\circle*{1}}          \put(133,22){$u_2$}
\put(60,5){\line(2,0){15}}      \put(60,5){\line(1,1){15}}
\put(60,5){\line(0,2){15}}      \put(60,20){\line(2,0){15}}
\put(75,5){\line(2,0){15}}      \put(75,5){\line(0,2){15}}
\put(75,20){\line(2,0){15}}     \put(82.5,12.5){\line(1,1){7.5}}
\put(82.5,12.5){\line(-1,1){7.5}} \put(82.5,12.5){\line(1,-1){7.5}}
\put(82.5,12.5){\line(-1,-1){7.5}}\put(105,12.5){\line(2,1){15}}
\put(105,12.5){\line(2,-1){15}}\put(112.5,12.5){\line(1,1){7.5}}
\put(112.5,12.5){\line(1,-1){7.5}}\put(120,5){\line(2,0){15}}
\put(120,5){\line(0,2){15}}\put(120,5){\line(1,1){15}}
\put(120,5){\line(1,2){7.5}}\put(120,20){\line(2,0){15}}
\put(135,5){\line(0,2){15}}
\end{picture}\\[2mm]

\unitlength=1mm
\begin{picture}(150,25)
\put(7.5,12.5){\circle*{1}}     \put(15,12.5){\circle*{1}}
\put(22.5,12.5){\circle*{1}}    \put(30,5){\circle*{1}}
\put(30,20){\circle*{1}}        \put(37.5,12.5){\circle*{1}}
\put(45,5){\circle*{1}}         \put(60,12.5){\circle*{1}}
\put(67.5,12.5){\circle*{1}}    \put(75,5){\circle*{1}}
\put(75,20){\circle*{1}}        \put(82.5,12.5){\circle*{1}}
\put(90,5){\circle*{1}}         \put(90,20){\circle*{1}}
\put(105,5){\circle*{1}}        \put(105,12.5){\circle*{1}}
\put(112.5,12.5){\circle*{1}}   \put(120,5){\circle*{1}}
\put(120,20){\circle*{1}}       \put(127.5,12.5){\circle*{1}}
\put(135,5){\circle*{1}}          \put(7.5,12.5){\line(3,1){22.5}}
\put(7.5,12.5){\line(3,-1){22.5}} \put(15,12.5){\line(2,1){15}}
\put(15,12.5){\line(2,-1){15}}    \put(22.5,12.5){\line(1,1){7.5}}
\put(22.5,12.5){\line(1,-1){7.5}} \put(30,5){\line(0,1){15}}
\put(30,5){\line(1,1){7.5}}       \put(30,20){\line(1,-1){7.5}}
\put(30,5){\line(1,0){15}}        \put(45,5){\line(-1,1){7.5}}
\put(60,12.5){\line(2,1){15}}     \put(60,12.5){\line(2,-1){15}}
\put(67.5,12.5){\line(1,1){7.5}}  \put(67.5,12.5){\line(1,-1){7.5}}
\put(75,5){\line(0,1){15}}        \put(75,5){\line(1,0){15}}
\put(82.5,12.5){\line(1,1){7.5}}  \put(82.5,12.5){\line(-1,-1){7.5}}
\put(82.5,12.5){\line(1,-1){7.5}} \put(82.5,12.5){\line(-1,1){7.5}}
\put(90,5){\line(0,1){15}}        \put(105,5){\line(0,1){7.5}}
\put(105,5){\line(1,0){15}}       \put(105,12.5){\line(2,1){15}}
\put(105,12.5){\line(2,-1){15}}   \put(112.5,12.5){\line(1,-1){7.5}}
\put(112.5,12.5){\line(1,1){7.5}} \put(120,5){\line(0,1){15}}
\put(120,5){\line(1,0){15}}       \put(120,5){\line(1,1){7.5}}
\put(127.5,12.5){\line(1,-1){7.5}}\put(127.5,12.5){\line(-1,1){7.5}}
\put(2.5,11.5){$u_7$}\put(17,11.5){$u_6$}\put(24.5,11.5){$u_5$}
\put(23,22){$u_2(x,u)$}\put(21,1){$u_1(y,v,w)$}\put(39.5,11.5){$u_3$}
\put(43,1){$u_4$}\put(55,11.5){$u_6$}\put(62.5,11.5){$u_5$}
\put(70,22){$u_2(u)$}\put(66,1){$u_1(y,v,w)$}\put(85.9,11.5){$u_3(x)$}
\put(88,22){$u_7$}\put(88,1){$u_4$}\put(100,11.5){$u_4$}
\put(100,1){$u_7$}\put(107.5,11.5){$u_6$}\put(111,1){$u_1(y,v,w)$}
\put(113,22){$u_3(x,u)$}\put(129.5,11.5){$u_2$}\put(133,1){$u_5$}
\put(28,-4){$T_4$}\put(73,-4){$T_5$}\put(118,-4){$T_6$}
\end{picture}\\

\unitlength=1mm
\begin{picture}(150,25)
\put(7.5,12.5){\circle*{1}}       \put(2.5,11.5){$u_5$}
\put(15,20){\circle*{1}}          \put(8,22){$u_2(x,w)$}
\put(15,5){\circle*{1}}           \put(13,1){$u_1$}
\put(30,20){\circle*{1}}          \put(25,22){$u_6(u)$}
\put(30,5){\circle*{1}}           \put(23,1){$u_3(y,v)$}
\put(45,5){\circle*{1}}           \put(43,1){$u_7$}
\put(45,20){\circle*{1}}          \put(43,22){$u_4$}
\put(60,12.5){\circle*{1}}        \put(55,11.5){$u_7$}
\put(67.5,12.5){\circle*{1}}      \put(62.5,11.5){$u_4$}
\put(75,5){\circle*{1}}           \put(68,1){$u_1(y,v)$}
\put(75,20){\circle*{1}}          \put(66,22){$u_3(x,u,w)$}
\put(82.5,12.5){\circle*{1}}      \put(84.5,11.5){$u_2$}
\put(90,5){\circle*{1}}           \put(88,1){$u_5$}
\put(90,20){\circle*{1}}          \put(88,22){$u_6$}
\put(105,5){\circle*{1}}          \put(102,1){$u_5$}
\put(105,20){\circle*{1}}         \put(102,22){$u_2$}
\put(112.5,20){\circle*{1}}       \put(107.5,22){$u_3(x,u)$}
\put(120,5){\circle*{1}}          \put(111,1){$u_1(y,v,w)$}
\put(127.5,20){\circle*{1}}       \put(125.5,22){$u_6$}
\put(135,20){\circle*{1}}         \put(133,22){$u_4$}
\put(135,5){\circle*{1}}          \put(133,1){$u_7$}
\put(7.5,12.5){\line(1,-1){7.5}}  \put(7.5,12.5){\line(1,1){7.5}}
\put(15,5){\line(1,0){15}}        \put(15,5){\line(1,1){15}}
\put(15,5){\line(0,1){15}}        \put(15,20){\line(1,0){15}}
\put(30,5){\line(1,0){15}}        \put(30,5){\line(1,1){15}}
\put(30,5){\line(0,1){15}}        \put(30,20){\line(1,0){15}}
\put(45,5){\line(0,1){15}}        \put(60,12.5){\line(2,-1){15}}
\put(60,12.5){\line(2,1){15}}     \put(67.5,12.5){\line(1,-1){7.5}}
\put(67.5,12.5){\line(1,1){7.5}}  \put(75,5){\line(1,0){15}}
\put(75,5){\line(0,1){15}}        \put(82.5,12.5){\line(1,-1){7.5}}
\put(82.5,12.5){\line(1,1){7.5}}  \put(82.5,12.5){\line(-1,-1){7.5}}
\put(82.5,12.5){\line(-1,1){7.5}} \put(75,20){\line(1,0){15}}
\put(105,5){\line(1,0){30}}       \put(105,5){\line(0,1){15}}
\put(105,20){\line(1,0){30}}      \put(120,5){\line(-1,1){15}}
\put(120,5){\line(-1,2){7.5}}     \put(120,5){\line(1,2){7.5}}
\put(120,5){\line(1,1){15}}       \put(135,5){\line(0,1){15}}
\put(28,-4){$T_7$}\put(73,-4){$T_8$}\put(118,-4){$T_9$}
\end{picture}\\

\unitlength=1mm
\begin{picture}(150,25)
\put(7.5,12.5){\circle*{1}}       \put(2.5,11.5){$u_5$}
\put(15,20){\circle*{1}}          \put(8,22){$u_1(x,u)$}
\put(15,5){\circle*{1}}           \put(13,1){$u_2$}
\put(30,20){\circle*{1}}          \put(28,22){$u_6$}
\put(30,5){\circle*{1}}           \put(21,1){$u_3(y,v,w)$}
\put(45,5){\circle*{1}}           \put(43,1){$u_7$}
\put(37.5,12.5){\circle*{1}}      \put(39.5,11.5){$u_4$}
\put(60,5){\circle*{1}}           \put(58,1){$u_6$}
\put(67.5,12.5){\circle*{1}}      \put(62.5,11.5){$u_4$}
\put(75,5){\circle*{1}}           \put(66,1){$u_1(y,v,w)$}
\put(75,20){\circle*{1}}          \put(73,22){$u_3$}
\put(82.5,12.5){\circle*{1}}      \put(85.9,11.5){$u_2(x,u)$}
\put(90,5){\circle*{1}}           \put(88,1){$u_5$}
\put(90,20){\circle*{1}}          \put(88,22){$u_7$}
\put(105,5){\circle*{1}}          \put(102,1){$u_6$}
\put(105,20){\circle*{1}}         \put(102,22){$u_5$}
\put(112.5,20){\circle*{1}}       \put(107.5,22){$u_2(x,u)$}
\put(120,5){\circle*{1}}          \put(111,1){$u_1(y,v,w)$}
\put(127.5,20){\circle*{1}}       \put(125.5,22){$u_3$}
\put(135,20){\circle*{1}}         \put(133,22){$u_4$}
\put(135,5){\circle*{1}}          \put(133,1){$u_7$}
\put(7.5,12.5){\line(1,-1){7.5}}  \put(7.5,12.5){\line(1,1){7.5}}
\put(15,5){\line(1,0){15}}        \put(15,5){\line(0,1){15}}
\put(15,20){\line(1,-1){15}}      \put(15,20){\line(1,0){15}}
\put(30,5){\line(1,0){15}}        \put(30,5){\line(0,1){15}}
\put(37.5,12.5){\line(-1,1){7.5}} \put(37.5,12.5){\line(-1,-1){7.5}}
\put(37.5,12.5){\line(1,-1){7.5}} \put(60,5){\line(1,0){15}}
\put(60,5){\line(1,1){7.5}}       \put(67.5,12.5){\line(1,-1){7.5}}
\put(67.5,12.5){\line(1,1){7.5}}  \put(75,5){\line(1,0){15}}
\put(75,5){\line(0,1){15}}        \put(82.5,12.5){\line(1,-1){7.5}}
\put(82.5,12.5){\line(1,1){7.5}}  \put(82.5,12.5){\line(-1,-1){7.5}}
\put(82.5,12.5){\line(-1,1){7.5}} \put(90,5){\line(0,1){15}}
\put(105,5){\line(1,0){30}}       \put(105,5){\line(1,2){7.5}}
\put(105,20){\line(1,0){30}}      \put(135,5){\line(-1,2){7.5}}
\put(120,5){\line(-1,1){15}}      \put(120,5){\line(-1,2){7.5}}
\put(120,5){\line(1,2){7.5}}      \put(120,5){\line(1,1){15}}
\put(28,-4){$T_{10}$}\put(73,-4){$T_{11}$}\put(118,-4){$T_{12}$}
\end{picture}
\begin{center}
Figure 1 (The all 2-trees on 7 vertices)
\end{center}

Let $T(k)=K_2\vee\overline{K_{k-2}}$, where $\overline{K_{k-2}}$ is
the complement of $K_{k-2}$ and $\vee$ denotes join operation.
Clearly, $T(k)$ is an 2-tree on $k$ vertices and has $k-2$ ears, and
every ear attaches to the unique edge of $K_2$. We recursively
define a new graph $F(k)$ as follows. Let $F(3)=K_3$, and let
$V(F(k-1))=\{x_1,\ldots,x_{k-1}\}$ for $k\ge 4$. Define $F(k)$ to be
the graph obtained from $F(k-1)$ by adding a new vertex $x_k$ and
joining $x_k$ to $x_{k-2}$ and $x_{k-1}$. Clearly, $F(k)$ is an
2-tree on $k$ vertices.

{\bf Lemma 3.1}\ \ {\it Let $G\in \{T_1,\ldots,T_{12}\}$. Then

 (1)\ \ there exist $x,y\in V(G)$ so that $G-\{x,y\}$ is a subgraph of $P_5$;

 (2)\ \ there exist $u,v\in V(G)$ so that
$G-\{u,v\}$ is a subgraph of $K_3\cup K_2$;

 (3)\ \ if $G\not=T_1$ (that is $T(7)$), then there exists $w\in
V(G)$ so that $G-w$ is a subgraph of $F(6)$.}

{\bf Proof.}\ \ By Figure 1, it is easy to check that $T_i-
\{x,y\}$ is a subgraph of $P_5$ for $1\le i\le 12$ and
$T_i-\{u,v\}$ is a subgraph of $K_3\cup K_2$ for $1\le i\le
12$. Moreover, $T_i-w$ is a subgraph of $F(6)$ for
$2\le i\le 12$.\ \ $\Box$

{\bf Lemma 3.2}\ \ {\it Let $n\ge30$ and $\pi=(d_1,\ldots,d_n)\in
GS_n$ with $\sigma(\pi)>6n-10$. Then $d_1\ge 6$, $d_3\ge 5$, $d_6\ge
4$ and $d_7\ge 3$.}

{\bf Proof.}\ \ If $d_1\le 5$, then $\sigma(\pi)\le 5n<6n-10$, a
contradiction. Hence $d_1\ge 6$. If $d_3\le4$, then
$\sigma(\pi)\leq2(n-1)+4(n-2)=6n-10$, a contradiction. Hence
$d_3\ge5$. If $d_6\le3$, by Theorem 2.2, then
$$\begin{array}{lll}
\sigma(\pi)&=&\sum\limits_{i=1}^{n}d_i
=\sum\limits_{i=1}^{5}d_i+\sum\limits_{i=6}^{n}d_i
\leq(5\times 4+\sum\limits_{i=6}^{n}\min\{5,d_i\})+\sum\limits_{i=6}^{n}d_i\\
&=&20+2\sum\limits_{i=6}^{n}d_i\leq 20+6(n-5)\\
&=&6n-10,\end{array}$$ a contradiction. Hence $d_6\ge4$. If
$d_7\le2$, by Theorem 2.2, then
$$\begin{array}{lll}
\sigma(\pi)&=&\sum\limits_{i=1}^{n}d_i= \sum\limits_{i=1}^{6}d_i+\sum\limits_{i=7}^{n}d_i
\leq(6\times 5+\sum\limits_{i=7}^{n}\min\{6,d_i\})+\sum\limits_{i=7}^{n}d_i\\
&=&30+2\sum\limits_{i=7}^{n}d_i\leq 30+4(n-6)\\
&<&6n-10,\end{array}$$ a contradiction. Hence $d_7\ge3$. \ \ $\Box$

We now define a new graph $G(7)$ as follows: Let
$V(K_{4})=\{v_1,v_2,v_3,v_4\}$ and $G(7)$ be the graph obtained from
$K_4$ by adding new vertices $x_1,x_2,x_3$, joining $x_i$ to
$v_1,\ldots,v_{i+1}$ for $1\le i\le 3$ and joining $x_1$ to $x_2$.

{\bf Lemma 3.3}\ \ {\it If $G$ is any 2-tree on 7 vertices, then
$G(7)$ contains $G$ as a subgraph.}

{\bf Proof.}\ \ Clearly, $T(7)$ is a subgraph of $G(7)$ and $F(6)$
is a subgraph of $G(7)-v_1$. Assume $G\not=T(7)$. By Lemma 3.1(3),
there exists $w\in V(G)$ so that $G-w$ is a subgraph of $F(6)$.
Hence $G-w$ is a subgraph of $G(7)-v_1$. We can see that $G(7)$
contains $G$ by putting $w$ on $v_1$. \ \ $\Box$

Let $k=7,n\ge30$ and $\pi=(d_1,\ldots,d_n)\in GS_n$ satisfy $d_2\ge6,
d_3\ge 5,d_6\ge4$ and $d_n \ge3$. We now define sequence $\pi_0,\pi_1,
\ldots,\pi_7$ as follows. Let $\pi_0=(d_1,\ldots,d_5,d_6-1,d_7-1,d_8,
\ldots,d_n)$. We define the sequence
$$\pi_1=(d_{2}^{(1)},\ldots,d_{7}^{(1)},d_{8}^{(1)},\ldots,d_{n}^{(1)})$$
from $\pi_0$ by deleting $d_1$, decreasing the first $d_1$ remaining
nonzero terms each by one unity, and then reordering the last $n-7$
terms to be non-increasing.

For $2\leq i\leq 7$ and $i\not=6$, we define the sequence
$$
\pi_i=(d_{i+1}^{(i)},\ldots,d_{7}^{(i)},d_{8}^{(i)},\ldots,d_{n}^{(i)})$$
from
$$\pi_{i-1}=(d_{i}^{(i-1)},\ldots,d_{7}^{(i-1)},d_{8}^{(i-1)},\ldots,
d_{n}^{(i-1)})$$ by deleting $d_i^{(i-1)}$, decreasing the first
$d_i^{(i-1)}$ remaining nonzero terms each by one unity, and then
reordering the last $n-7$ terms to be non-increasing.

For $i=6$ , we define the sequence
$$
\pi_6=(d_{7}^{(6)},d_{8}^{(6)},\ldots,d_{n}^{(6)})$$ from
$$\pi_{5}=(d_{6}^{(5)},d_{7}^{(5)},d_{8}^{(5)},\ldots,d_{n}^{(5)})$$
by deleting $d_6^{(5)}$, decreasing the first $d_6^{(5)}$ remaining
nonzero terms each by one unity except for the term $d_7^{(5)}$, and
then reordering
 the last $n-7$ terms to be non-increasing.

By the definition of $\pi_7$, the following Proposition 3.1 is obvious.

{\bf Proposition 3.1}\ \ {\it Let $k=7, n\ge 30$ and
$\pi=(d_1,\ldots,d_n)\in GS_n$ satisfy $d_2\ge6, d_3\ge 5,d_6\ge 4$
and $d_n \ge 3$. If $\pi_7$ is graphic, then $\pi$ has a realization
containing $G(7)$.}

{\bf Lemma 3.4}\ \ {\it Let $k=7$, $n\ge 30$ and
$\pi=(d_1,\ldots,d_7,d_8,\ldots,d_n)\in GS_n$ satisfy $d_2\ge6,
d_3\ge 5, d_6\ge4$ and $d_n \ge3$. Let
$\pi_1'=(d_1',\ldots,d_{n-1}')$ be the residual sequence obtained
from $\pi$ by laying off $d_1$, and let
$\rho=(\rho_1,\ldots,\rho_{n-2})$ be the residual sequence obtained
from $\pi_1'$ by laying off the term $d_2-1$. If $\pi$ satisfies one
of (a)--(c), where

(a)\ \ $d_1=d_2=n-1$,

(b)\ \  $d_1=n-1$, $d_2\leq n-2$ and $d_7>d_{d_2+2}$,

(c)\ \ $d_1\leq n-2$, $d_7>d_{d_2+2}$ and $d_7-d_{d_1+2}\geq 2$,\\
then $\rho_1=d_3-2,\rho_2=d_4-2,\ldots,\rho_5=d_7-2$.}

{\bf Proof.}\ \ If $\pi$ satisfies (a), then
$\rho=(d_3-2,\ldots,d_n-2)$, and so
$\rho_1=d_3-2,\rho_2=d_4-2,\ldots,\rho_5=d_7-2$.

If $\pi$ satisfies (b), then $\pi_1'=(d_2-1,d_3-1,\ldots,d_n-1)$. By
$d_7-2\ge d_{d_2+2}-1$, we have
$\rho_1=d_3-2,\rho_2=d_4-2,\ldots,\rho_5=d_7-2$.

Assume that $\pi$ satisfies (c). If $d_{d_2+2}>d_{d_1+2}$, then
$d_{d_2+2}-1\ge d_{d_1+2}$, and hence
$d_1'=d_2-1,\ldots,d_{d_2+1}'=d_{d_2+2}-1$. By $d_7>d_{d_2+2}$,
we have $d_7-2\ge d_{d_2+2}-1$, implying that
$\rho_1=d_3-2,\rho_2=d_4-2,\ldots,\rho_5=d_7-2$. If
$d_{d_2+2}=\cdots=d_{d_1+2}$, then $d_{d_2+2}-1<d_{d_1+2}$. By
$d_7-d_{d_1+2}\geq 2$, we have
$d_1'=d_2-1,\ldots,d_6'=d_7-1$ and $d_{d_2+1}'\le
d_{d_1+2}$, implying that
$\rho_1=d_3-2,\rho_2=d_4-2,\ldots,\rho_5=d_7-2$.\ \ $\Box$

{\bf Lemma 3.5}\ \ {\it Let $k=7$, $n\geq 30$ and
$\pi=(d_1,\ldots,d_n)
 \in GS_n$ satisfy $d_2\ge6, d_3\ge 5,d_6\ge4$ and $d_n \ge3$. For
each
$\pi_i=(d_{i+1}^{(i)},\ldots,d_{7}^{(i)},d_{8}^{(i)},\ldots,d_{n}^{(i)})$,
let $s_i=\max\{j|d_{8}^{(i)}-d_{7+j}^{(i)}\leq1\}$.

(1)\ \  If $\pi$ satisfies (d) or (e), where (d) $d_1\leq n-2$,
$d_{7}>d_{d_2+2}$ and $d_{7}-d_{d_1+2}\leq 1$ and (e) $d_1\leq n-2$,
$d_{7}=d_{d_2+2}$ and $d_{d_2+2}=d_{d_1+2}$, then
$d_{7+r}^{(7)}=d_{7+r}$ for $r>s_{7}$.

(2)\ \ If $\pi$ satisfies (f) or (g), where (f) $d_1=n-1$, $d_2\leq
n-2$ and $d_{7}=d_{d_2+2}$ and (g) $d_1\leq n-2$, $d_{7}=d_{d_2+2}$
and $d_{d_2+2}>d_{d_1+2}$, then $d_{7+r}^{(7)}=d_{7+r}^{(1)}$ for
$r>s_{7}$.}

{\bf Proof.}\ \ (1) If $\pi$ satisfies (d) or (e), then $7+s_0\ge
d_1+2$. Since $d_{8}^{(i-1)}-d_{3t+s_{i-1}}^{(i-1)}\le 1$ implies
that $d_{8}^{(i)}-d_{7+s_{i-1}}^{(i)}\leq1$ for $1\leq i\leq 7$, we
have that $s_7\geq s_6\geq\cdots\geq s_0\geq d_1+2-7$. By
$\min\{d_{8}^{(i-1)}-1,\ldots,d_{d_i+1}^{(i-1)}-1,d_{d_i+2}^{(i-1)},\ldots,
d_{7+s_{i-1}}^{(i-1)}\}\geq d_{8}^{(i-1)}-2\geq
d_{7+s_{i-1}+1}^{(i-1)}\geq\cdots\geq d_{n}^{(i-1)}$, we have that
$d_{7+s_{i-1}+m}^{(i)}=d_{7+s_{i-1}+m}^{(i-1)}$ for $m\geq1$. Thus,
$d_{7+r}^{(i)}=d_{7+r}^{(i-1)}$ for $r>s_{i}$. This implies that
$d_{7+r}^{(7)}=d_{7+r}$ for $r>s_7$.

(2) If $\pi$ satisfies (f) or (g), then $s_7 \geq s_6\geq
\cdots \geq s_1\ge s_0\ge d_2+2-7$. Since
$\min\{d_{8}^{(i-1)}-1,\ldots,d_{d_i+1}^{(i-1)}-1,d_{d_i+2}^{(i-1)},\ldots,
d_{7+s_{i-1}}^{(i-1)}\}\geq d_{7+1}^{(i-1)}-2\geq
d_{7+s_{i-1}+1}^{(i-1)}\geq\cdots\geq d_{n}^{(i-1)}$ for $i\geq 2$,
we have that $d_{7+s_{i-1}+m}^{(i)}=d_{7+s_{i-1}+m}^{(i-1)}$ for
$i\geq 2$ and $m\geq1$. Thus, $d_{7+r}^{(i)}=d_{7+r}^{(i-1)}$ for
$i\geq 2$ and $r>s_{i}$. This implies that
$d_{7+r}^{(7)}=d_{7+r}^{(1)}$ for $r>s_7$.\ \ $\Box$

{\bf Lemma 3.6}\ \ {\it Let $n\ge 6$ and $\pi=(d_1,\ldots,d_n)\in
GS_n$ with $d_n\ge 1$ and $\sigma(\pi)>2n$. Then $\pi$ has a
realization containing $K_3\cup K_2$ or $P_5$.}

{\bf Proof.}\ \ To the contrary, we assume that $\pi$ has no
realization containing $K_3\cup K_2$ or $P_5$. By $n\ge
6,\sigma(\pi)>2n$ and Theorem 2.7(1), $\pi$ has a realization $G$
containing $K_3$. Let $V(G)=\{v_1,\ldots,v_n\}$ so that the subgraph
induced by $\{v_1,v_2,v_3\}$ is $K_3$. Then
$E(G-\{v_1,v_2,v_3\})=\emptyset$ as $G$ contains no $K_3\cup K_2$.
If there are $i,j\in \{4,\ldots,n\}$ so that $v_i$ has two neighbors
in $\{v_1,v_2,v_3\}$ and $v_j$ has one neighbor in
$\{v_1,v_2,v_3\}$, then $G$ contains $P_5$, a contradiction. Hence
$v_i$ has only one neighbor in $\{v_1,v_2,v_3\}$ for each $i\in
\{4,\ldots,n\}$. Thus $|E(G)|=n-3+3=n$, implying $\sigma(\pi)=2n$, a
contradiction.\ \ $\Box$

 If $\pi=(d_1,\ldots,d_n)\in GS_n$ has a realization containing every
2-tree on $k$ vertices, then $\pi$ is {\it potentially} $A'(k)$-{\it
graphic}. If $\pi$ has a realization in which the subgraph induced
by the $k$ vertices of highest degree contains every 2-tree on $k$
vertices, then $\pi$ is {\it potentially} $A''(k)$-{\it graphic}.
Clearly, if $\pi$ is potentially $A''(k)$-graphic, then $\pi$ is
potentially $A'(k)$-graphic.

{\bf Lemma 3.7}\ \ {\it Let $n\ge 30$ and $\pi=(d_1,\ldots,d_n)\in
GS_n$ with  $d_2\ge6$ and $d_n\ge3$. If $\sigma(\pi)>6n-10$, then
$\pi$ is potentially $A''(7)$-graphic.}

{\bf Proof.}\ \ Let $\pi_1'=(d_1',\ldots,d_{n-1}')$ be the residual
sequence obtained from $\pi$ by laying off $d_1$, and let
$\rho=(\rho_1,\ldots,\rho_{n-2})$ be the residual sequence obtained
from $\pi_1'$ by laying off the term $d_2-1$. Then
$n-2>6,\rho_{n-2}\ge 3-2=1$ and
$\sigma(\rho)=\sigma(\pi)-2d_1-2d_2+2>6n-10-4(n-1)+2=2(n-2)$. By
Lemma 3.6 and Theorem 2.4, $\rho$ has a realization $G_1$ in which
the subgraph induced by the vertices with degrees
$\rho_1,\ldots,\rho_5$ contains $K_3\cup K_2$ or $P_5$. Denote $F$
to be the subgraph induced by the vertices with degrees
$\rho_1,\ldots,\rho_5$ in $G_1$, and let $F'$ be the graph obtained
from $F$ by adding two vertices $x,y$ such that $x,y$ are adjacent
to each vertex of $F$. Since $F$ contains $K_3\cup K_2$ or $P_5$, by
Lemma 3.1 (1) and (2), we can see that $F'$ contains every 2-tree on
7 vertices.

If $\pi$ satisfies one of (a)--(c), by Lemma 3.4, then we have
$\rho_1=d_3-2,\rho_2=d_4-2,\ldots,\rho_5=d_7-2$. Now by the
definitions of $\rho$ and $\pi_1'$, it is easy to get that $\pi$ has
a realization $G'$ in which the subgraph induced by the vertices
with degrees $d_1,\ldots,d_7$ contains $F'$. In other words, $\pi$
is potentially $A''(7)$-graphic.

We assume that $\pi$ satisfies one of (d)--(g). If $d_7\ge 11$, then
by Theorem 2.5, $\pi$ has a realization containing $K_7$, and hence
$\pi$
 is potentially $A''(7)$-graphic by Theorem 2.4. Assume that $d_7\le 10$.
  By Lemma 3.2, we have $d_3\ge 5$ and $d_6 \ge 4$.
 It is enough to prove that $\pi_7$ is graphic by Theorem 2.4, Lemma 3.3
 and Proposition 3.1. If $\pi$ satisfies (d) or (e), by Lemma 3.5(1), then
$$\pi_7=(d_{8}^{(7)},\ldots,d_{7+s_7}^{(7)},d_{7+s_{7}+1},\ldots,d_n).$$
If $\pi$ satisfies (f) or (g), by Lemma 3.5(2), then
$$\pi_7=(d_{8}^{(7)},\ldots,d_{7+s_7}^{(7)},
d_{7+s_{7}+1}^{(1)},\ldots,d_n^{(1)}).$$ If $s_7<n-7$,  by
$d_{8}^{(7)}\le d_7\le 10$ and $d_n\ge d_n^{(1)}\ge d_n-1\ge 2$,
then we have that
$$
\frac{1}{2}\left\lfloor\frac{(10+2+1)^2}{4}\right\rfloor\le
\frac{13^2}{8}<22<n-7.$$ By Theorem 2.3, $\pi_7$ is graphic. If
$s_7=n-7$, then $d_{8}^{(7)}-d_n^{(7)}\le 1$. Denote $d_n^{(7)}=m$.
If $m=0$, then by $d_8^{(7)}\le 1$ and $\sigma(\pi_7)$ being even,
$\pi_7$ is clearly graphic. If $m\ge 1$, then $d_{8}^{(7)}\le m+1$,
and hence
$$\frac{1}{m}\left\lfloor\frac{(m+1+m+1)^2}{4}\right\rfloor=
\frac{(m+1)^2}{m}\leq m+3\le 10+3<n-7.$$ By Theorem 2.3, $\pi_7$ is
also graphic.\ \ $\Box$

{\bf Lemma 3.8}\ \ {\it Let $n=30+s$ with $0\le s\le 63$, and let
$\pi=(d_1,\ldots,d_n)\in GS_n$ with $\sigma(\pi)>6n-6+126-2s$. Then
$\pi$ is potentially $A'(7)$-graphic.}

{\bf Proof.}\ \ We use induction on $s$. If $s=0$, then $n=30$ and
$\sigma(\pi)\ge 6\times 30-6+126-0+2=300+2=2\times 30\times
(7-2)+2$. By Theorem 2.6, $\pi$ has a realization containing $K_7$,
and hence $\pi$ is potentially $A'(7)$-graphic. Suppose now that
$1\le s\le 63$. Then $\sigma(\pi)>6n-6$. If $d_2\le 5$, then
$\sigma(\pi)\le (n-1)+5(n-1)=6n-6$, a contradiction. Hence $d_2\ge
6$. If $d_n\ge 3$, then $\pi$ is potentially $A'(7)$-graphic by
Lemma 3.7. If $d_n\le 2$, then the residual sequence
$\pi_{n}'=(d_{1}',\ldots,d_{n-1}')$ obtained by laying off $d_n$
from $\pi$ satisfies
$\sigma(\pi_{n}')=\sigma(\pi)-2d_n>6n-6+126-2s-4=6(n-1)-6+126-2(s-1)$.
By the induction hypothesis, $\pi_{n}'$ is potentially
$A'(7)$-graphic, and hence so is $\pi$.\  \ $\Box$

{\bf Proof of Theorem 1.4.}\ \ Let $k\in\{3,4,5,7\}$, $n\ge N(k)$
and $\pi=(d_1,\ldots,d_n)\in GS_n$ with $\sigma(\pi)>(k-1)(n-1)$. We
first prove that $\pi$ is potentially $A'(k)$-graphic. If $k=3$, by
$\sigma(\pi)>(k-1)(n-1)=2n-2$, $n\ge N(3)=6$ and Theorem 2.7(1),
then $\pi$ has a realization containing $K_3$, and hence $\pi$ is
potentially $A'(3)$-graphic. If $k=4$,  by $\sigma(\pi)
>(k-1)(n-1)=3n-3$, $n\ge N(4)=7$ and Theorem 2.7(2), then $\pi$ has a
realization containing $K_4-e$. This implies that $\pi$ is
potentially $A'(4)$-graphic. If $k=5$, by
$\sigma(\pi)>(k-1)(n-1)=4n-4$, $n\ge N(5)=24$ and Theorem 2.7(3),
then $\pi$ has a realization containing $K_5-E(P_3)$. Since
$K_5-E(P_3)$ contains every 2-tree on 5 vertices ($F(5)$ and $T(5)$
are the only two 2-trees on 5 vertices), we have that $\pi$ is
potentially $A'(5)$-graphic. Assume that $k=7$. Then
$\sigma(\pi)>6n-6$ and $d_2\ge 6$. We use induction on $n\ge
N(7)=93$. If $n=93$, by Lemma 3.8 ($s=63$), then $\pi$ is
potentially $A'(7)$-graphic. Assume that $n\ge 94$. If $d_n\ge 3$,
then $\pi$ is potentially $A''(7)$-graphic by Lemma 3.7. If $d_n\le
2$, then the residual sequence $\pi_{n}'$ satisfies
$\sigma(\pi_{n}')=\sigma(\pi)-2d_n>6n-6-4>6(n-1)-6$. By the
induction hypothesis, $\pi_{n}'$ is potentially $A'(7)$-graphic, and
hence so is $\pi$.

We now show that $(k-1)(n-1)$ is the best possible lower bound in
Theorem 1.4. For $k\in\{3,4,5,7\}$, let $\pi=(n-1,(k-2)^{n-1})$,
where the symbol $x^y$ stands for $y$ consecutive terms $x$.
Clearly, $\pi$ is graphic and $\sigma(\pi)=(k-1)(n-1)$. Since $T(k)$
has two vertices with degree $k-1$ and $d_2=k-2$ in $\pi$, we have
that every realization of $\pi$ contains no $T(k)$ as a subgraph.
Thus, we must require that $\sigma(\pi)>(k-1)(n-1)$ in Theorem 1.4.\
\ $\Box$

\section*{4. Proof of Theorem 1.5}

\hskip\parindent In order to prove Theorem 1.5, we also need the
following lemmas.

{\bf Lemma 4.1}\ \ {\it Let $G$ be any 2-tree on $k\ge 6$ vertices and
$G\not=T(k)$. Let $B(G)$ be the set of all ears in $G$ and
$C(G)=\{e(u)|u\in B(G)\}$. Then $|C(G)|\ge 2$.}

{\bf Proof.}\ \ If $|C(G)|=1$, let $C(G)=\{xy\}$, then $u$ attaches
to $xy$ for each $u\in B(G)$. Let $G'=G-B(G)$. Since $G\not=F(k)$,
we have that $|V(G')|\ge 3$, $G'$ is an 2-tree and each vertex of
$V(G')-\{x,y\}$ has degree at least 3 in $G'$. This implies that
$G'\not=K_3$, and $x$ and $y$ are exactly two ears in $G'$ by
Theorem 2.8(1). This is impossible by Theorem 2.8(3).\ \ $\Box$

{\bf Lemma 4.2}\ \ {\it Let $G$ be any 2-tree on $k\ge 6$ vertices.
Let $xy\in C(G)$ so that $xy$ is attached by $s$ ears
$z,z_1,\ldots,z_{s-1}$. Then $G-\{x,y,z\}$ is a spanning subgraph of
some 2-tree on $k-3$ vertices.}

{\bf Proof.}\ \ Assume that $G\not=T(k)$. Let $G'=G-
\{z,z_1,\ldots,z_{s-1}\}$. Then $G'$ is an 2-tree on $k-s$ vertices.
By $|C(G)|\ge 2$ (Lemma 4.1), we have that $k-s\ge4$. By Theorem
2.8(6), $G'$ can be constructed from $xy$ by repeatedly adding a new
vertex and making it adjacent to the two ends of an edge in the
graph formed so far. In the process of constructing $G'$ from $xy$,
we let $y'$ be the first vertex that is attached to $xy$. Since $xy$
can not be attached by an ear in $G'$, we have that $d_{G'}(y')\ge
3$. This implies that $xy'$ or $yy'$ must be attached by a new
vertex. Let $x'$ be the first vertex that is attached to $xy'$ or
$yy'$. Without loss of generality, we assume that $x'$ is attached
to $xy'$. Let $\{x_1,\ldots,x_t\}$ be the subset of $V(G')$ so that
$x_i$ is attached to $xx'$ for $i=1,\ldots,t$ and
$\{y_1,\ldots,y_{t'}\}$ be the subset of $V(G')$ so that $y_j$ is
attached to $yy'$ for $j=1,\ldots,t'$. Denote
$$G''=G'-\{xx_1,\ldots,xx_t\}-\{yy_1,\ldots,yy_{t'}\}+\{y'x_1,\ldots,
y'x_t\}+\{x'y_1,\ldots,x'y_{t'}\}-\{xy,xy'\}.$$ In $G''$, we
identify the vertex $x$ to the vertex $x'$ and the vertex $y$ to the
vertex $y'$, and attach $z_r$  to $x'y'$ for $r=1,\ldots,s-1$, the
resulting graph is denoted by $G'''$. Then $G'''$ is an 2-tree on
$k-3$ vertices. Clearly, $G-\{x,y,z\}$ is a spanning subgraph of
$G'''$.

If $G=T(k)$, then $G-\{x,y,z\}$ is an independent set on $k-3$
vertices, and hence $G-\{x,y,z\}$ is a spanning subgraph of some
2-tree on $k-3$ vertices.\ \ $\Box$

The following Lemma 4.3 is a further result of Lemma 4.2 for $k\ge
8$. For $X\subseteq V(G)$ and $v\in V(G)$, the neighborhood of $v$
in $X$ is denoted by $N_X(v)$.

{\bf Lemma 4.3}\ \ {\it Let $G$ be any 2-tree on $k\ge 8$ vertices.
Let $xy\in C(G)$ so that $xy$ is attached by $s$ ears
$z,z_1,\ldots,z_{s-1}$. Then $G-\{x,y,z\}$ is a spanning subgraph of
some 2-tree $T$ on $k-3$ vertices and $T\not=T(k-3)$.}

{\bf Proof.}\ \ Denote $H=G-\{x,y,z\}$. If $s\ge 2$, then $G-z_1$ is
an 2-tree on $k-1$ vertices, and by Lemma 4.2, $(G-z_1)-\{x,y,z\}$
is a spanning subgraph of some 2-tree $T'$ on $k-4$ vertices. We can
get a new 2-tree $T$ by attaching a new vertex $z_1$ to a suitable
edge of $T'$ so that $T\not=T(k-3)$. Clearly, $H$ is a spanning
subgraph of $T$.

Assume that $xy$ is attached by only an ear for each $xy\in C(G)$.
To the contrary, we assume that the 2-tree on $k-3$ vertices that
contains $H$ as a spanning subgraph is the only $T(k-3)$. Let
$V(T(k-3))=\{u,v,w_1,\ldots,w_{k-5}\}$ so that $uv\in C(T(k-3))$ and
$uv$ is attached by $k-5$ ears $w_1,\ldots,w_{k-5}$. If $w_pu\notin
E(H)$, where $1\le p\le k-5$, then $H$ is a spanning subgraph of an
2-tree $T=T(k-3)-w_pu+w_pw_q$, where $1\le q\le k-5$ and $q\not=p$.
This is impossible as $T\not=T(k-3)$. Hence  $w_pu\in E(H)$ for
$1\le p\le k-5$. Similarly, $w_pv\in E(H)$ for $1\le p\le k-5$. Thus
$H=T(k-3)$ if $uv\in E(H)$ or $H=K_{2,k-5}$ if $uv\notin E(H)$.

If $H=T(k-3)$, by $|E(G)|=2k-3$ and $|E(H)|=2(k-3)-3=2k-9$, the
number of edges in $G$ joining $\{x,y\}$ to $V(H)$ is
$e_G(\{x,y\},V(H))=(2k-3)-(2k-9)-3=3$. If $|N_H(x)|=3$ or
$|N_H(y)|=3$, then $G$ is not 2-connected, a contradiction by
Theorem 2.8(5). Hence $|N_H(x)|\le 2$ and $|N_H(y)|\le 2$. Denote
$W=\{w_1,\ldots,w_{k-5}\}$. If $|N_W(x)|=2$ or $|N_W(y)|=2$, without
loss of generality, we assume that $|N_W(x)|=2$ and let
$w_ix,w_jx\in E(G)$. Then $xw_iuw_jx$ is a chordless cycle of length
4 in $G$, a contradiction by Theorem 2.8(4). If
$|N_W(x)|=|N_W(y)|=1$, let $w_rx,w_ty\in E(G)$, then $xw_ruw_tyx$ is
a cycle of length 5 in $G$ with at most one chord by
$e_G(\{x,y\},V(H))=3$. This implies that $G$ contains a chordless
cycle of length at least 4, a contradiction. If
$|N_W(x)|+|N_W(y)|\le 1$, then $uv\in C(G)$ and $uv$ is attached by
at least $k-5-1\ge 2$ ears, a contradiction. If $H=K_{2,k-5}$, then
$H$ contains a chorless cycle of length 4, and hence so is $G$, a
contradiction.  \ \ $\Box$

For $\pi=(d_1,\ldots,d_n) \in GS_n$, we always assume that
$\pi_1'=(d_1',\ldots,d_{n-1}')$ is the residual sequence obtained
from $\pi$ by laying off $d_1$ and $\rho=(\rho_1,\ldots,\rho_{n-2})$
is the residual sequence obtained from $\pi_1'$ by laying off the
term $d_2-1$. The proof of the following Lemma 4.4 is similar to
that of Lemma 3.4, we omit it here.

{\bf Lemma 4.4}\ \ {\it Let $k\ge6$, $n\geq 6k$ and
$\pi=(d_1,\ldots,d_n) \in GS_n$. If $\pi$ satisfies one of (a)--(c),
where

(a)\ \ $d_1=d_2=n-1$,

(b)\ \  $d_1=n-1$, $d_2\leq n-2$ and $d_k>d_{d_2+2}$,

(c)\ \ $d_1\leq n-2$, $d_k>d_{d_2+2}$ and $d_k-d_{d_1+2}\geq 2$,\\
then $\rho_1=d_3-2,\rho_2=d_4-2,\ldots,\rho_{k-2}=d_k-2$.}

Theorem 1.3 implies that the lower bound
$2\lfloor\frac{2k}{3}\rfloor
n-2n-\lfloor\frac{2k}{3}\rfloor^2+\lfloor\frac{2k}{3}\rfloor+1-(-1)^i$
in Theorem 1.5 is the best possible. Thus, the proof of Theorem 1.5
can be divided into the proofs of the following three theorems in
terms of the value of $k=3t$, $3t+1$ and $3t+2$.

{\bf Theorem 4.1} \ \ {\it If $t\ge2$, $n\ge20t^2-t$ and
$\pi=(d_1,\ldots,d_n)\in GS_n$ with $\sigma(\pi)>4tn-2n-4t^2+2t$,
then $\pi$ has a realization containing every 2-tree on $3t$
vertices.}

{\bf Theorem 4.2} \ \ {\it If $t\ge3$, $n\ge20t^2+23t+5$ and
$\pi=(d_1,\ldots,d_n)\in GS_n$ with $\sigma(\pi)>4tn-2n-4t^2+2t+2$,
then $\pi$ has a realization containing every 2-tree on $3t+1$
vertices.}

{\bf Theorem 4.3} \ \ {\it If $t\ge2$, $n\ge20t^2+31t+12$ and
$\pi=(d_1,\ldots,d_n)\in GS_n$ with $\sigma(\pi)>4tn-4t^2-2t$, then
$\pi$ has a realization containing every 2-tree on $3t+2$ vertices.}

\section*{4.1. Proof of Theorem 4.1}

\hskip\parindent {\bf Lemma 4.1.1}\ \ {\it Let $t\ge2,n\ge18t$ and
$\pi=(d_1,\ldots,d_n)\in GS_n$ with $\sigma(\pi)>4tn-2n-4t^2+2t$.
Then

 (1) $d_i\ge 3t-\lceil\frac{i}{2}\rceil$ for
$i=1,\ldots,2t$;

 (2) $d_i\ge2(3t+1-i)$ for $i=2t+1,\ldots,3t$.}

{\bf Proof.}\ \ (1) If there is an even $s$ with $2\le s\le 2t$ such
that $d_s\le 3t-\lceil\frac{s}{2}\rceil-1=3t-\frac{s}{2}-1$, then
$$\begin{array}{lll}
\sigma(\pi)&\leq&(s-1)(n-1)+(3t- \frac{s}{2}-1)(n-s+1)\\
&=&\frac{s^2}{2}-s(3t-\frac{n}{2}+\frac{1}{2})+3tn-2n+3t.
\end{array}$$
Denote $f(s)=\frac{s^2}{2}-s(3t-\frac{n}{2}+\frac{1}{2})+3tn-2n+3t$.
Since $2\le s\le 2t$, we have that
$$\begin{array}{lll}
\sigma(\pi)&\leq&f(s) \leq \max\{f(2),f(2t)\}\\
&=&\max\{3tn-n-3t+1,4tn-2n-4t^2+2t\}\\
&=&\max\{4tn-2n-4t^2+2t-(n-4t+1)(t-1),4tn-2n-4t^2+2t\}\\
&=& 4tn-2n-4t^2+2t,
\end{array}$$ a contradiction.

If there is an odd $s$ with $1\le s\le 2t-1$ such that $d_s\le
3t-\lceil\frac{s}{2}\rceil-1=k-\frac{s+1}{2}-1$, then
$$\begin{array}{lll}
\sigma(\pi)&\leq&(s-1)(n-1)+(3t- \frac{s+1}{2}-1)(n-s+1)\\
&=&\frac{s^2}{2}-s(3t-\frac{n}{2})+3tn+3t-\frac{5n}{2}-\frac{1}{2}.
\end{array}$$
Denote
$g(s)=\frac{s^2}{2}-s(3t-\frac{n}{2})+3tn+3t-\frac{5n}{2}-\frac{1}{2}$.
Since $1\leq s \leq 2t-1$, we have that
$$\begin{array}{lll}
\sigma(\pi)&\leq&g(s) \leq \max\{g(1),g(2t-1)\}\\
&=&\max\{3tn-2n-9t^2+9t-2,4tn-3n-4t^2+4t\}\\
&=&\max\{4tn-2n-4t^2+2t-[t(n-2)+5t(t-1)+2],4tn-2n-4t^2+2t-(n-2t)\}\\
&\le&4tn-2n-4t^2+2t,
\end{array}$$ a contradiction.

(2) If there is $s$ with $2t+1\le s\le 3t$ such that $d_s\le
6t-2s+1$,  by Theorem 2.2, then
$$\begin{array}{lll}
\sigma(\pi)&=&\sum\limits_{i=1}^{n}d_i
=\sum\limits_{i=1}^{s-1}d_i+\sum\limits_{i=s}^{n}d_i
\leq((s-2)(s-1)+\sum\limits_{i=s}^{n}\min\{s-1,d_i\})+\sum\limits_{i=s}^{n}d_i\\
&=&(s-2)(s-1)+2\sum\limits_{i=s}^{n}d_i\leq(s-2)(s-1)+2(6t-2s+1)(n-s+1)\\
&=&5s^2-(12t+4n+9)s+12tn+12t+2n+4.
\end{array}$$
Denote $h(s)=5s^2-(12t+4n+9)s+12tn+12t+2n+4$. Since $2t+1\le s\le
3t$, we have that
$$\begin{array}{lll}
\sigma(\pi)&\leq&h(s) \leq \max\{h(2t+1),h(3t)\}\\
&=&\max\{4tn-2n-4t^2+2t,2n+9t^2-15t+4\}\\
&=&\max\{4tn-2n-4t^2+2t,4tn-2n-4t^2+2t-[(n-4t)(4t-4)+3t^2+t-4]\}\\
&=&4tn-2n-4t^2+2t,
\end{array}$$ a contradiction.\ \ $\Box$

 We now define a new graph $G(3t)$ as follows: Let $V(K_{2t})=\{v_1,v_2,\ldots,v_{2t}
 \}$ and $G(3t)$ be the graph obtained from $K_{2t}$ by adding new
 vertices $x_1,x_2,\ldots,x_t$ and joining $x_i$ to $v_1,v_2,\ldots,v_{2i}$
 for $1\leq i \leq t$.

{\bf Lemma 4.1.2}\ \ {\it If $G$ is any 2-tree on $3t$ vertices,
then $G(3t)$ contains $G$.}

{\bf Proof.}\ \  We use induction on $t$. It is obvious for $t=1$.
Assume $t\ge 2$. Let $xy\in C(G)$ so that $xy$ is attached by the
ear $z$. Denote $H=G-\{x,y,z\}$. By Lemma 4.2, $H$ is a spanning
subgraph of some 2-tree $G'$ on $3(t-1)$ vertices. Let
$$M=G(3t)-\{v_1,v_2\}-\{x_1\}.$$
Then $M=G(3(t-1))$. By the induction hypothesis, $G(3(t-1))$
contains $G'$. This implies that $G(3(t-1))$ contains $H$. Putting
$x,y$ and $z$ on $v_1,v_2$ and $x_1$ respectively, we can see that
$G(3t)$ contains $G$.\ \ $\Box$

 Let $t\ge2$, $n\geq 18t$ and
$\pi=(d_1,\ldots,d_n)
 \in GS_n$ satisfy

 (i)\ \ $d_i\ge 3t-\lceil\frac{i}{2}\rceil$ for
$i=1,\ldots,2t$,

(ii)\ \ $d_{2t+1}\ge2t$, and

(iii)\ \ $d_n\ge2t-1$.

We now define sequence $\pi_0,\pi_1,\ldots,\pi_{3t}$ as follows. Let
$\pi_0=\pi$. We define the sequence
$$
\pi_1=(d_{2}^{(1)},\ldots,d_{3t}^{(1)},d_{3t+1}^{(1)},\ldots,d_{n}^{(1)})$$
from $\pi_0$ by deleting $d_1$, decreasing the first $d_1$ remaining
nonzero terms each by one unity, and then reordering the last $n-3t$
terms to be non-increasing.

For $2\leq i\leq 3t$, we define the sequence
$$
\pi_i=(d_{i+1}^{(i)},\ldots,d_{3t}^{(i)},d_{3t+1}^{(i)},\ldots,d_{n}^{(i)})$$
from
$$\pi_{i-1}=(d_{i}^{(i-1)},\ldots,d_{3t}^{(i-1)},d_{3t+1}^{(i-1)},\ldots,
d_{n}^{(i-1)})$$ by deleting $d_i^{(i-1)}$, decreasing the first
$d_i^{(i-1)}$ remaining nonzero terms each by one unity, and then
reordering the last $n-3t$ terms to be non-increasing.

By the definition of $\pi_{3t}$, the following Proposition 4.1.1 is obvious.

{\bf Proposition 4.1.1}\ \ {\it Let $t\ge2$, $n\geq 18t$ and
$\pi=(d_1,\ldots,d_n)\in GS_n$ satisfy (i)--(iii). If $\pi_{3t} $ is
graphic, then $\pi$ has a realization containing $G(3t)$.}

{\bf Lemma 4.1.3}\ \ {\it Let $t\ge2$, $n\geq 18t$ and
$\pi=(d_1,\ldots,d_n)
 \in GS_n$ satisfy (i)--(iii). For each
$\pi_i=(d_{i+1}^{(i)},\ldots,d_{3t}^{(i)},d_{3t+1}^{(i)},\ldots,d_{n}^{(i)})$,
let $s_i=\max\{j|d_{3t+1}^{(i)}-d_{3t+j}^{(i)}\leq1\}$.

(1)\ \ If $\pi$ satisfies (d) or (e), where (d) $d_1\leq n-2$,
$d_{3t}>d_{d_2+2}$ and $d_{3t}-d_{d_1+2}\leq 1$ and (e) $d_1\leq
n-2$, $d_{3t}=d_{d_2+2}$ and $d_{d_2+2}=d_{d_1+2}$, then
$d_{3t+r}^{(3t)}=d_{3t+r}$ for $r>s_{3t}$.

(2)\ \ If $\pi$ satisfies (f) or (g), where (f) $d_1=n-1$, $d_2\le
n-2$ and $d_{3t}=d_{d_2+2}$ and (g) $d_1\le n-2$, $d_{3t}=d_{d_2+2}$
and $d_{d_2+2}>d_{d_1+2}$, then $d_{3t+r}^{(3t)}=d_{3t+r}^{(1)}$ for
$r>s_{3t}$.}

{\bf Proof.}\ \  The proof of Lemma 4.1.3 is similar to that of
Lemma 3.5.\ \ $\Box$

{\bf Lemma 4.1.4}\ \ {\it Let $t\ge 1$, $n\ge 18t$ and
$\pi=(d_1,\ldots,d_n)\in GS_n$ with $d_n\ge 2t-1$ and
$\sigma(\pi)>4tn-2n-4t^2+2t$. Then $\pi$ is potentially
$A''(3t)$-graphic.}

{\bf Proof.}\ \ We use induction on $t$.  If $t=1$, then
$\sigma(\pi)>2n-2$. By $n\ge 18$ and Theorem 2.7(1), $\pi$ has a
realization containing $K_3$, implying that $\pi$ is potentially
$A''(3)$-graphic by Theorem 2.4. Assume $t\ge 2$. By the definition
of $\rho=(\rho_1,\ldots,\rho_{n-2})$, we can see that $n-2\ge
18(t-1)$, $\rho_{n-2}\ge (2t-1)-2
 =2(t-1)-1$ and $\sigma(\rho)=\sigma(\pi)-2d_1
 -2d_2+2>4tn-2n-4t^2+2t-4(n-1)+2=4(t-1)(n-2)-2(n-2)-4(t-1)^2+2(t-1)$.
 By the induction hypothesis, $\rho$ has a realization $G_1$ in which the
 subgraph $F$ induced by the vertices with degrees $\rho_1,\ldots,\rho_{3(t-1)}$
 contains every 2-tree on $3(t-1)$ vertices. Denote $F'$ to be the graph obtained from $F$ by adding
 three new vertices $x,y,u$ such that $x,y$ are adjacent to each
 vertex of $F$ and $xy,xu,yu\in E(F')$.

{\bf Claim}\ \ {\it $F'$ contains every 2-tree on $3t$ vertices.}

{\it Proof of Claim.}\ \ Let $G$ be any 2-tree on $3t$ vertices.
Take $xy\in C(G)$ and $u\in B(xy)$, and denote $H=G-\{x,y,u\}$. By
Lemma 4.2, it is easy to get that $H$ is a spanning subgraph of some
2-tree on $3(t-1)$ vertices. Since $F$ contains every 2-tree on
$3(t-1)$ vertices, we have that $F$ contains $H$. By the definition
of $F'$, we can see that $F'$ contains $G$. By the arbitrary of $G$,
$F'$ contains every 2-tree on $3t$ vertices. This proves Claim.\ \
$\Box$

If $\pi$ satisfies one of (a)--(c), by Lemma 4.4 (the case of
$k=3t$), then
$\rho_1=d_3-2,\rho_2=d_4-2,\ldots,\rho_{3t-2}=d_{3t}-2$. This is
implies that $\pi$ has a realization $G'$ in which the subgraph
induced by the vertices with degrees $d_1,\ldots,d_{3t}$ contains
$F'$. Thus by Claim, $\pi$ is potentially $A''(3t)$-graphic.

 We now assume that $\pi$ satisfies one of (d)--(g). If
 $d_{3t}\ge 6t-3$, by Theorem 2.5, then $\pi$ has a realization containing
 $K_{3t}$, and hence $\pi$ is potentially $A''(3t)$-graphic by Theorem 2.4.
 Assume that $d_{3t}\leq6t-4$. By Lemma 4.1.1,
 we have $d_i\ge 3t-\lceil \frac{i}{2} \rceil$ for $i=1,\ldots,2t$ and
 $d_{2t+1}\geq2t$. It is enough to prove that $\pi_{3t}$ is graphic by
 Theorem 2.4, Lemma 4.1.2 and Proposition 4.1.1. If $\pi$ satisfies (d) or (e), by
 Lemma 4.1.3(1), then
$$\pi_{3t}=(d_{3t+1}^{(3t)},\ldots,d_{3t+s_{3t}}^{(3t)},d_{3t+s_{3t}+1},\ldots,d_n).$$
If $\pi$ satisfies (f) or (g), by Lemma 4.1.3(2), then
$$\pi_{3t}=(d_{3t+1}^{(3t)},\ldots,d_{3t+s_{3t}}^{(3t)},
d_{3t+s_{3t}+1}^{(1)},\ldots,d_n^{(1)}).$$ If $s_{3t}<n-3t$, by
$d_{3t+1}^{(3t)}\le d_{3t}\leq 6t-4$ and $d_n\ge d_n^{(1)}\ge
d_n-1\ge 2t-2\ge 2$,  then we have that
$$\begin{array}{lll}
\frac{1}{2t-2}\left\lfloor\frac{(6t-4+2t-2+1)^2}{4}\right\rfloor&\le&\frac{(8t-5)^2}{8(t-1)}\\
&=&\frac{64t^2-80t+25}{8(t-1)}\\
&=&\frac{64t(t-1)-16(t-1)+9}{8(t-1)}\\
&\leq&8t\leq n-3t.
\end{array}$$
By Theorem 2.3, $\pi_{3t}$ is graphic. If $s_{3t}=n-3t$, then
$d_{3t+1}^{(3t)}-d_n^{(3t)}\le 1$. Denote $d_n^{(3t)}=m$. If $m=0$,
by $d_{3t+1}^{(3t)}\le 1$ and $\sigma(\pi_{3t})$ being even, then
$\pi_{3t}$ is clearly graphic. If $m\ge 1$, then $d_{3t+1}^{(3t)}\le
m+1$, and hence
$$\frac{1}{m}\left\lfloor\frac{(m+1+m+1)^2}{4}\right\rfloor=\frac{(m+1)^2}{m}\leq m+3\le 6t-4+3\le
n-3t.$$ By Theorem 2.3, $\pi_{3t}$ is also graphic.\ \ $\Box$

{\bf Lemma 4.1.5}\ \ {\it Let $t\ge 2$, $n=18t+s$ with $0\le s\le 20t^2-19t$ and let
$\pi=(d_1,\ldots,d_n)\in GS_n$ with
$\sigma(\pi)>4tn-2n+36t^2-36t-2s$. Then
$\pi$ is potentially $A'(3t)$-graphic.}

{\bf Proof.}\ \ We use induction on $s$. If $s=0$, then $n=18t$ and
$\sigma(\pi)\ge 4tn-2n+36t^2-36t+2=2n(3t-2)+2$. By Theorem 2.6,
$\pi$ has a realization containing $K_{3t}$, and hence $\pi$ is
potentially $A'(3t)$-graphic. Suppose now that $1\le s\le
20t^2-19t$. Then $\sigma(\pi)>4tn-2n-4t^2+2t$. If $d_n\ge 2t-1$,
then $\pi$ is potentially $A'(3t)$-graphic by Lemma 4.1.4. If
$d_n\le2t-2$, then $\pi_{n}'=(d_{1}',\ldots,d_{n-1}')$ satisfies
$\sigma(\pi_{n}')=\sigma(\pi)-2d_n>4tn-2n+36t^2-36t-2s-2(2t-2)=
4t(n-1)-2(n-1)+36t^2-36t-2(s-1)$. By the induction hypothesis,
$\pi_{n}'$ is potentially $A'(3t)$-graphic, and hence so is $\pi$.\ \ $\Box$

{\bf Proof of Theorem 4.1.}\ \ Let $t\ge2,n\ge20t^2-t$ and
$\pi=(d_1,\ldots,d_n)\in GS_n$ with $\sigma(\pi)>4tn-2n-4t^2+2t$. We
only need to prove that $\pi$ is potentially $A'(3t)$-graphic. We
use induction on $n$. If $n=20t^2-t$, by Lemma 4.1.5
$(s=20t^2-19t)$, then $\pi$ is potentially $A'(3t)$-graphic. Assume
that $n\ge20t^2-t+1$. If $d_n\ge 2t-1$, by Lemma 4.1.4, then $\pi$
is potentially $A'(3t)$-graphic. If $d_n\le 2t-2$, then $\pi_{n}'$
satisfies
$\sigma(\pi_{n}')=\sigma(\pi)-2d_n\ge4tn-2n-4t^2+2t-2(2t-2)>4t(n-1)-2(n-1)-4t^2+2t$.
By the induction hypothesis, $\pi_{n}'$ is potentially
$A'(3t)$-graphic, and hence so is $\pi$.\ \ $\Box$

\section*{4.2. Proof of Theorem 4.2}

\hskip\parindent {\bf Lemma 4.2.1}\ \ {\it Let $t\ge 3, n\ge 18t+6$ and
$\pi=(d_1,\ldots,d_n)\in GS_n$ with $\sigma(\pi)>4tn-2n-4t^2+2t+2$.
Then

(1) $d_i\ge 3t+1-\lceil\frac{i}{2}\rceil$ for $i=1,\ldots,2t-3$ and
$d_{2t-1}\ge 2t+1$;

(2) $d_i\ge 2(3t+2-i)$ for $i=2t+2,\ldots,3t+1$.}

{\bf Proof.}\ \ (1) If there is an even $s$ with $2\le s\le 2t-4$
such that $d_s\le (3t+1)-\lceil\frac{s}{2}\rceil-1=3t-\frac{s}{2}$,
then
$$\begin{array}{lll}
\sigma(\pi)&\leq&(s-1)(n-1)+(3t-\frac{s}{2})(n-s+1)\\
&=&\frac{s^2}{2}-s(3t-\frac{n}{2}+\frac{3}{2})+3tn-n+3t+1.
\end{array}$$
Denote
$f(s)=\frac{s^2}{2}-s(3t-\frac{n}{2}+\frac{3}{2})+3tn-n+3t+1$.
Since $2\le s\le 2t-4$, we have that
$$\begin{array}{lll}
\sigma(\pi)&\leq&f(s) \leq \max\{f(2),f(2t-4)\}\\
&=&\max\{3tn-3t,4tn-3n-4t^2+4t+15\}\\
&=&\max\{4tn-2n-4t^2+2t+2-[(n-4t)(t-3)+n-7t+2],\\
& & 4tn-2n-4t^2+2t+2-(n-2t-13)\}\\
&\le&4tn-2n-4t^2+2t+2,
\end{array}$$ a contradiction. Hence $d_i\ge 3t+1-\lceil\frac{i}{2}\rceil$
for  even $s$ with $2\le s\le 2t-4$.

If there is an odd $s$ with $1\le s\le2t-1$ such that $d_s\le
(3t+1)-\lceil\frac{s}{2}\rceil-1=3t-\frac{s+1}{2}$, then
$$\begin{array}{lll}
\sigma(\pi)&\leq&(s-1)(n-1)+(3t-\frac{s+1}{2})(n-s+1)\\
&=&\frac{s^2}{2}-s(3t-\frac{n}{2}+1)+3tn+3t-\frac{3n}{2}+\frac{1}{2}.
\end{array}$$
Denote
$g(s)=\frac{s^2}{2}-s(3t-\frac{n}{2}+1)+3tn+3t-\frac{3n}{2}+\frac{1}{2}$.
Since $1\leq s \leq 2t-1$, we have that
$$\begin{array}{lll}
\sigma(\pi)&\leq&g(s) \leq \max\{g(1),g(2t-1)\}\\
&=&\max\{3tn-n,4tn-2n-4t^2+2t+2\}\\
&=&\max\{4tn-2n-4t^2+2t+2-(n-4t-2)(t-1),4tn-2n-4t^2+2t+2\}\\
&=&4tn-2n-4t^2+2t+2,
\end{array}$$ a contradiction. Hence $d_i\ge 3t+1-\lceil\frac{i}{2}\rceil$
for odd $s$ with $1\le s\le 2t-1$, that is, $d_i\ge
3t+1-\lceil\frac{i}{2}\rceil$ for odd $s$ with $1\le s\le 2t-3$ and
$d_{2t-1}\ge 3t+1-\lceil\frac{2t-1}{2}\rceil=2t+1$.

(2) If there is an $s$ with $2t+2\le s\le 3t+1$ such that $d_s\le
6t-2s+3$, by Theorem 2.2, then
$$\begin{array}{lll}
\sigma(\pi)
&=&\sum\limits_{i=1}^{n}d_i
=\sum\limits_{i=1}^{s-1}d_i+\sum\limits_{i=s}^{n}d_i
\leq((s-2)(s-1)+\sum\limits_{i=s}^{n}\min\{s-1,d_i\})+\sum\limits_{i=s}^{n}d_i\\
&=&(s-2)(s-1)+2\sum\limits_{i=s}^{n}d_i\leq(s-2)(s-1)+2(6t-2s+3)(n-s+1)\\
&=&5s^2-(12t+4n+13)s+12tn+12t+6n+8.
\end{array}$$
Denote $h(s)=5s^2-(12t+4n+13)s+12tn+12t+6n+8$. Since $2t+2\le s\le
3t+1$, we have that
$$\begin{array}{lll}
\sigma(\pi)&\leq&h(s) \leq \max\{h(2t+2),h(3t+1)\}\\
&=&\max\{4tn-2n-4t^2+2t+2,2n+9t^2-9t\}\\
&=&\max\{4tn-2n-4t^2+2t+2,4tn-2n-4t^2+2t-[4(n-4t)(t-1)+3t^2-5t]\}\\
&=&4tn-2n-4t^2+2t+2,
\end{array}$$ a contradiction.\ \ $\Box$

We now define a new graph $G(3t+1)$ for $t\ge3$ as follows. Let
$V(K_{2t+1})=\{v_1,\ldots,v_{2t+1}\}$ and $G(3t+1)$ be the graph
obtained from the graph $K_{2t+1}-v_{2t-2}v_{2t}$ by adding new
vertices $x_1,\ldots,x_{t}$ and joining $x_i$ to $v_1,\ldots,v_{2i}$
for $1\le i\le t$. \\[2mm]
\unitlength=1mm
\begin{picture}(150,40)
\put(60,5){\circle*{1}}            \put(55,4){$u_2$}
\put(60,27.5){\circle*{1}}         \put(55,29.5){$u_3$}
\put(75,12.5){\circle*{1}}         \put(70,11.5){$u_1$}
\put(75,35){\circle*{1}}           \put(70,37){$u_4$}
\put(105,5){\circle*{1}}            \put(107,4){$u_5$}
\put(105,20){\circle*{1}}           \put(107,19){$u_6$}
\put(105,35){\circle*{1}}           \put(107,37){$u_7$}
\put(60,5){\line(3,0){45}}         \put(60,5){\line(3,1){45}}
\put(60,5){\line(3,2){45}}         \put(60,5){\line(2,1){15}}
\put(60,5){\line(0,1){22.5}}       \put(75,12.5){\line(-1,1){15}}
\put(75,12.5){\line(0,1){22.5}}    \put(60,27.5){\line(2,1){15}}
\put(75,12.5){\line(4,1){30}}      \put(75,12.5){\line(4,3){30}}
\put(75,12.5){\line(4,-1){30}}     \put(75,35){\line(1,0){30}}
\put(75,35){\line(2,-1){30}}       \put(60,27.5){\line(6,-1){45}}
\put(60,27.5){\line(6,1){45}}
\end{picture}
\begin{center}
Figure 2 (The graph $M$)
\end{center}

{\bf Lemma 4.2.2}\ \ {\it If $t\ge 3$ and $G$ is any 2-tree on
$3t+1$ vertices, then $G(3t+1)$ contains $G$.}

{\bf Proof.}\ \ We use induction on $t$. Let $xy\in C(G)$ so that
$xy$ is attached by the ear $z$. Denote $H=G-\{x,y,z\}$. By Lemma
4.3, $H$ is a spanning subgraph of some 2-tree $T$ on $3(t-1)+1$
vertices and $T\not=T(3(t-1)+1)$. Assume $t=3$. By Figure 1--2, we
can see that the graph $M$ in Figure 2 contains any 2-tree $T$ on
$7$ vertices, where $T\not=T(7)$. Since $G(10)-\{v_1,v_2,x_1\}$ is
isomorphic to $M$, we have that  $H$ is a spanning subgraph of
$G(10)-\{v_1,v_2,x_1\}$. Putting $x,y$ and $z$ on $v_1,v_2$ and
$x_1$ respectively, we can see that $G(10)$ contains $G$ as a
subgraph. If $t\ge4$, then $G(3t+1) -\{v_1,v_2,x_1\}=G(3(t-1)+1)$
contains $T$, and hence contains $H$. Putting $x,y$ and $z$ on
$v_1,v_2$ and $x_1$ respectively, we can see that $G(3t+1)$ contains
$G$.\ \ $\Box$

Let $t\ge3,\ n\ge 18t+6$ and $\pi=(d_1,\ldots,d_n)\in GS_n$ satisfy

 (i)\ \  $d_i\ge 3t+1-\lceil\frac{i}{2}\rceil$ for $i=1,\ldots,2t-3$ and
$d_{2t-1}\ge 2t+1$,

(ii)\ \ $d_{2t+2}\ge 2t$, and

(iii)\ \ $d_n\ge 2t-1$.

We now define sequence $\pi_0,\pi_1,\ldots,\pi_{3t+1}$ as follows.
Let $\pi_0=\pi$. We define the sequence
$$\pi_1=(d_{2}^{(1)},\ldots,d_{3t+1}^{(1)},d_{3t+2}^{(1)},\ldots,d_{n}^{(1)})$$
from $\pi_0$ by deleting $d_1$, decreasing the first $d_1$ remaining
nonzero terms each by one unity, and then reordering the last
$n-(3t+1)$ terms to be non-increasing.

For $2\le i\le 3t+1$ and $i\not=2t-2$, we define the sequence
$$
\pi_i=(d_{i+1}^{(i)},\ldots,d_{3t+1}^{(i)},d_{3t+2}^{(i)},\ldots,d_{n}^{(i)})$$
from
$$\pi_{i-1}=(d_{i}^{(i-1)},\ldots,d_{3t+1}^{(i-1)},d_{3t+2}^{(i-1)},\ldots,d_{n}^{(i-1)})$$
by deleting $d_i^{(i-1)}$, decreasing the first
$d_i^{(i-1)}$ remaining nonzero terms each by one unity, and then
reordering the last $n-(3t+1)$ terms to be non-increasing.

For $i=2t-2$ and $d_{2t-2}\ge 2t+2$, the definition of $\pi_{2t-2}$
is as above. For $i=2t-2$ and $d_{2t-2}=2t+1$, we define the
sequence $$
\pi_{2t-2}=(d_{2t-1}^{(2t-2)},\ldots,d_{3t+1}^{(2t-2)},d_{3t+2}^{(2t-2)},\ldots,d_{n}^{(2t-2)})$$
from
$$\pi_{2t-3}=(d_{2t-2}^{(2t-3)},\ldots,d_{3t+1}^{(2t-3)},d_{3t+2}^{(2t-3)},\ldots,d_{n}^{(2t-3)})$$
by deleting $d_{2t-2}^{(2t-3)}$, decreasing the first
$d_{2t-2}^{(2t-3)}$ remaining nonzero terms each by one unity except
for the term $d_{2t}^{(2t-3)}$, and then reordering the last
$n-(3t+1)$ terms to be non-increasing.

By the definition of $\pi_{3t+1}$, the following Proposition 4.2.1 is obvious.

{\bf Proposition 4.2.1}\ \ {\it Let $t\ge3$, $n\geq 18t+6$ and
$\pi=(d_1,\ldots,d_n)\in GS_n$ satisfy (i)--(iii). If $\pi_{3t+1}$
is graphic, then $\pi$ has a realization containing $G(3t+1)$ on
those vertices with degrees $d_1,\ldots,d_{3t+1}$.}

{\bf Lemma 4.2.3}\ \ {\it Let $t\ge2$, $n\geq 18t+6$ and
$\pi=(d_1,\ldots,d_n)
 \in GS_n$ satisfy (i)--(iii). For each
$\pi_i=(d_{i+1}^{(i)},\ldots,d_{3t+1}^{(i)},d_{3t+2}^{(i)},\ldots,d_{n}^{(i)})$,
let $s_i=\max\{j|d_{3t+2}^{(i)}-d_{3t+1+j}^{(i)}\leq1\}$.

(1)\ \ If $\pi$ satisfies (d) or (e), where (d) $d_1\leq n-2$,
$d_{3t+1}>d_{d_2+2}$ and $d_{3t+1}-d_{d_1+2}\leq 1$ and (e) $d_1\leq
n-2$, $d_{3t+1}=d_{d_2+2}$ and $d_{d_2+2}=d_{d_1+2}$, then
$d_{3t+1+r}^{(3t+1)}=d_{3t+1+r}$ for $r>s_{3t+1}$.

(2)\ \ If $\pi$ satisfies (f) or (g), where (f) $d_1=n-1$, $d_2\leq
n-2$ and $d_{3t+1}=d_{d_2+2}$ and (g) $d_1\leq n-2$,
$d_{3t+1}=d_{d_2+2}$ and $d_{d_2+2}>d_{d_1+2}$, then
$d_{3t+1+r}^{(3t+1)}=d_{3t+1+r}^{(1)}$ for $r>s_{3t+1}$.}

{\bf Proof.}\ \ The proof of Lemma 4.2.3 is similar to that of Lemma
3.5.\ \ $\Box$

{\bf Lemma 4.2.4}\ \ {\it Let $n\ge 30$ and $\pi=(d_1,\ldots,d_n)\in
GS_n$ with $d_n\ge 3$ and $\sigma(\pi)>6n-10$. Then $\pi$ has a
realization $G'$ in which the subgraph induced by the vertices with
degrees $d_1,\ldots,d_7$ contains every 2-tree on 7 vertices except
for $T(7)$.}

{\bf Proof.}\ \ If $d_2\ge 6$, by Lemma 3.7, then $\pi$ is
potentially $A''(7)$-graphic, implying that $\pi$ has a realization
containing every 2-tree on 7 vertices on those vertices with degrees
$d_1,\ldots,d_7$. Assume that $d_2\le 5$. It follows from
$\sigma(\pi)>6n-10$ and $\sigma(\pi)$ being even that
$\pi=(n-1,5^{n-1})$ or $(n-2,5^{n-2},4)$ or $(n-1,5^{n-3},4^{2})$ or
$(n-1,5^{n-2},3)$ or $(n-3,5^{n-1})$. Let $G$ be any 2-tree on 7
vertices with $G\not=T(7)$. By Lemma 3.1(3), there exists $w\in
V(G)$ such that $G-w$ is a subgraph of $F(6)$. By the definition of
$F(n)$, we have that $d_{F(n)}(x_1)=d_{F(n)}(x_n)=2$,
$d_{F(n)}(x_2)=d_{F(n)}(x_{n-1})=3$ and
$d_{F(n)}(x_3)=\cdots=d_{F(n)}(x_{n-2})=4$, that is, the degree
sequence of $F(n)$ is $(4^{n-4} ,3^{2},2^{2})$.

If $\pi=(n-1,5^{n-1})$, then $\pi_1'=(4^{n-1})$. Clearly,
$G_1=F(n-1)+\{x_1x_{n-2},x_1x_{n-1},x_2x_{n-1}\}$ is a realization
of $\pi_1'$, and $G_1$ contains $F(6)$ on vertices $x_1,\ldots,x_6$.
Let $G'$ be the graph obtained from $G_1$ by adding a new vertex $x$
that is adjacent to each vertex of $G_1$. Clearly, $G'$ is a
realization of $\pi$, and $G'$ contains every 2-tree on 7 vertices
except for $T(7)$ on vertices $x,x_1,\ldots,x_6$ with degrees
$n-1,5,\ldots,5$.

If $\pi=(n-2,5^{n-2},4)$, then $\pi_1'=(4^{n-1})$. Clearly,
$G_1=F(n-1)+\{x_1x_{n-2},x_1x_{n-1},x_2x_{n-1}\}$ is a realization
of $\pi_1'$, and $G_1$ contains $F(6)$ on vertices $x_1,\ldots,x_6$.
Let $G'$ be the graph obtained from $G_1$ by adding a new vertex $x$
that is adjacent to $x_1,\ldots,x_{n-2}$. Clearly, $G'$ is a
realization of $\pi$, and $G'$ contains every 2-tree on 7 vertices
except for $T(7)$ on vertices $x,x_1,\ldots,x_6$ with degrees
$n-2,5,\ldots,5$.

If $\pi=(n-1,5^{n-3},4^{2})$, then $\pi_1'=(4^{n-3},3^{2})$.
Clearly, $G_1=F(n-1)+\{x_1x_{n-2},x_2x_{n-1}\}$ is a realization of
$\pi_1'$, and $G_1$ contains $F(6)$ on vertices $x_2,\ldots,x_7$.
Let $G'$ be the graph obtained from $G_1$ by adding a new vertex $x$
that is adjacent to each vertex of $G_1$. Clearly, $G'$ is a
realization of $\pi$, and $G'$ contains every 2-tree on 7 vertices
except for $T(7)$ on vertices $x,x_2,\ldots,x_7$ with degrees
$n-1,5,\ldots,5$.

If $\pi=(n-1,5^{n-2},3)$, then $\pi_1'=(4^{n-2},2)$. Clearly,
$$G_1=F(n-1)+\{x_1x_{n-2},x_1x_{n-3},x_2x_{n-2}\}-x_{n-2}x_{n-3}$$ is
a realization of $\pi_1'$, and $G_1$ contains $F(6)$ on vertices
$x_1,\ldots,x_6$. Let $G'$ be the graph obtained from $G_1$ by
adding a new vertex $x$ that is adjacent to each vertex of $G_1$.
Clearly, $G'$ is a realization of $\pi$, and $G'$ contains every
2-tree on 7 vertices except for $T(7)$ on vertices
$x,x_1,\ldots,x_6$ with degrees $n-1,5,\ldots,5$.

If $\pi=(n-3,5^{n-1})$, then $\pi_1'=(5^2,4^{n-3})$. Clearly,
$$G_1=F(n-1)+\{x_1x_{n-2},x_1x_{n-1},x_2x_{n-2},x_2x_{n-1}\}$$ is a
realization of $\pi_1'$, and $G_1$ contains $F(6)$ on vertices
$x_3,\ldots,x_8$. Let $G'$ be the graph obtained from $G_1$ by
adding a new vertex $x$ that is adjacent to
$x_1,x_3,x_4,\ldots,x_{n-3},x_{n-1}$. Clearly, $G'$ is a realization
of $\pi$, and $G'$ contains every 2-tree on 7 vertices except for
$T(7)$ on vertices $x,x_3,\ldots,x_8$ with degrees
$n-3,5,\ldots,5$.\ \ $\Box$

{\bf Lemma 4.2.5}\ \ {\it Let $t\ge 3$, $n\ge 18t+6$ and
$\pi=(d_1,\ldots,d_n)\in GS_n$ with $d_n\ge 2t-1$ and
$\sigma(\pi)>4tn-2n-4t^2+2t+2$. Then $\pi$ is potentially
$A''(3t+1)$-graphic.}

{\bf Proof.}\ \ We only need to prove the following two claims.

{\bf Claim 1}\ \ {\it If $\pi$ satisfies one of (d)--(g), then $\pi$
is potentially $A''(3t+1)$-graphic.}

{\it Proof of Claim 1.}\ \ If $d_{3t+1}\ge 6t-1$, by Theorem 2.5,
then $\pi$ has a realization containing $K_{3t+1}$, and hence $\pi$
is potentially $A''(3t+1)$-graphic by Theorem 2.4. Assume that
$d_{3t+1}\le 6t-2$. By Lemma 4.2.1, we have $d_i\ge
3t+1-\lceil\frac{i}{2}\rceil$ for $i=1,\ldots,2t-3$, $d_{2t-1}\ge
2t+1$ and $d_{2t+2}\ge 2t$. It is enough to prove that $\pi_{3t+1}$
is graphic by Lemma 4.2.2 and Proposition 4.2.1. If $\pi$ satisfies
(d) or (e), by Lemma 4.2.3(1), then
$$\pi_{3t+1}=(d_{3t+2}^{(3t+1)},\ldots,d_{3t+1+s_{3t+1}}^{(3t+1)},d_{3t+1+s_{3t+1}+1},\ldots,d_n).$$
If $\pi$ satisfies (f) or (g), by Lemma 4.2.3(2), then
$$\pi_{3t+1}=(d_{3t+2}^{(3t+1)},\ldots,d_{3t+1+s_{3t+1}}^{(3t+1)},d_{3t+1+s_{3t+1}+1}^{(1)},\ldots,d_{n}^{(1)}).$$
If $s_{3t+1}<n-(3t+1)$, by $d_{3t+2}^{(3t+1)}\le d_{3t+1}\le 6t-2$
and $d_n\ge d_n^{(1)}\ge d_n-1\ge 2t-2\ge 2$, then we have that
$$\begin{array}{lll}
\frac{1}{2t-2}\left\lfloor\frac{(6t-2+2t-2+1)^2}{4}\right\rfloor&\le &\frac{(6t-2+2t-2+1)^2}{4(2t-2)}\\
&=&\frac{(8t-3)^2}{8(t-1)}\\
&=&\frac{64t^2-48t+9}{8(t-1)}\\
&=&\frac{64t(t-1)+16(t-1)+25}{8(t-1)}\\
&\le&8t+2+2\le n-(3t+1).
\end{array}$$
By Theorem 2.3, $\pi_{3t+1}$ is graphic. If $s_{3t+1}=n-(3t+1)$,
then $d_{3t+2}^{(3t+1)}-d_n^{(3t+1)}\le 1$. Denote $d_n^{(3t+1)}=m$.
If $m=0$, then by $d_{3t+2}^{(3t+1)}\le 1$ and $\sigma(\pi_{3t+1})$
being even, $\pi_{3t+1}$ is clearly graphic. If $m\ge 1$, then
$d_{3t+2}^{(3t+1)}\le m+1$, and hence
$$\frac{1}{m}\left\lfloor\frac{(m+1+m+1)^2}{4}\right\rfloor=\frac{(m+1)^2}{m}\leq m+3\le 6t-2+3\le
n-(3t+1).$$ By Theorem 2.3, $\pi_{3t+1}$ is also graphic. This
proves Claim 1.\ \ $\Box$

{\bf Claim 2}\ \ {\it If $\pi$ satisfies one of (a)--(c), then $\pi$
is potentially $A''(3t+1)$-graphic.}

{\it Proof of Claim 2.}\ \ We use induction on $t$. Let $G$ be any
2-tree on $3t+1$ vertices, and let $x,y,z\in V(G)$ so that $xy\in
C(G)$ and $xy$ is attached by the ear $z$. Denote $H=G-\{x,y,z\}$.
By Lemma 4.3, $H$ is a spanning subgraph of some 2-tree $T$ on
$3(t-1)+1$ vertices and $T\not=T(3(t-1)+1)$.

If $t=3$, by the definition of $\rho$, we can see that $n-2\ge 58$,
$\rho_{n-2}\ge 5-2=3$ and
$\sigma(\rho)=\sigma(\pi)-2d_1-2d_2+2>10n-28-4(n-1)+2=6(n-2)-10$. By
Lemma 4.2.4, $\rho$ has a realization $G_1$ in which the subgraph
$F$ induced by the vertices with degrees $\rho_1,\ldots,\rho_7$
contains every 2-tree on 7 vertices except for $T(7)$. Denote $F'$
to be the graph obtained from $F$ by adding three new vertices
$x,y,z$ such that $x,y$ are adjacent to each vertex of $F$ and
$xy,xz,yz\in E(F')$. Since $F$ contains $T$ (and hence $H$), we have
that $F'$ contains $G$. By the arbitrary of $G$, $F'$ contains every
2-tree on 10 vertices. Moreover, by Lemma 4.4 (the case of $k=10$),
we have $\rho_1=d_3-2,\rho_2=d_4-2,\ldots,\rho_8=d_{10}-2$. This is
implies that $\pi$ has a realization $G'$ in which the subgraph
induced by the vertices with degrees $d_1,\ldots,d_{10}$ contains
$F'$. In other words, $\pi$ is potentially $A''(10)$-graphic.

If $t\ge 4$, by the definition of $\rho$, then we have that $n-2\ge
18(t-1)+6$, $\rho_{n-2}\ge 2t-1-2=2(t-1)-1$ and
$\sigma(\rho)=\sigma(\pi)-2d_1-2d_2+2>4tn-2n-4t^2+2t+2-4(n-1)+2
=4(t-1)(n-2)-2(n-2)-4(t-1)^2+2(t-1)+2$. If $\rho$ satisfies one of
(d)--(g),  by Claim 1, then $\rho$ is potentially
$A''(3(t-1)+1)$-graphic. If $\rho$ satisfies one of (a)--(c), by the
induction hypothesis, then $\rho$ is also potentially
$A''(3(t-1)+1)$-graphic. This implies that $\rho$ has a realization
$G_1$ in which the subgraph $F$ induced by the vertices with degrees
$\rho_1,\ldots,\rho_{3(t-1)+1}$ contains every 2-tree on $3(t-1)+1$
vertices. Denote $F'$ to be the graph obtained from $F$ by adding
three new vertices $x,y,z$ such that $x,y$ are adjacent to each
vertex of $F$ and $xy,xz,yz\in E(F')$. Since $F$ contains $T$ (and
hence $H$), we have that $F'$ contains $G$. By the arbitrary of $G$,
$F'$ contains every 2-tree on $3t+1$ vertices. Moreover, by Lemma
4.4 (the case of $k=3t+1$), we have
$\rho_1=d_3-2,\rho_2=d_4-2,\ldots,\rho_{3t-1}=d_{3t+1}-2$. This is
implies that $\pi$ has a realization $G'$ in which the subgraph
induced by the vertices with degrees $d_1,\ldots,d_{3t+1}$ contains
$F'$. In other words, $\pi$ is potentially $A''(3t+1)$-graphic. This
proves Claim 2.\ \ $\Box$

{\bf Lemma 4.2.6}\ \ {\it Let $t\ge 3$, $n=18t+6+s$ with $0\le s\le
20t^2+5t-1$ and $\pi=(d_1,\ldots,d_n)\in GS_n$ with $\sigma(\pi)
>4tn-2n+36t^2+12t-2s$. Then $\pi$ is potentially $A'(3t+1)$-graphic.}

{\bf Proof.}\ \ We use induction on $s$. If $s=0$, then $n=18t+6$
and $\sigma(\pi)\ge 4tn-2n+36t^2+12t+2=2n(3t+1-2)+2$. By Theorem
2.6, $\pi$ has a realization containing $K_{3t+1}$, and hence $\pi$
is potentially $A'(3t+1)$-graphic. Suppose now that $1\le s\le
20t^2+5t-1$. Then $\sigma(\pi)>4tn-2n-4t^2+2t+2$. If $d_n\ge 2t-1$,
then $\pi$ is potentially $A'(3t+1)$-graphic by Lemma 4.2.5. If
$d_n\le 2t-2$, then $\pi_{n}'=(d_{1}', \ldots,d_{n-1}')$ satisfies
$\sigma(\pi_{n}')=\sigma(\pi)-2d_n>4tn-2n+36t^2+
12t-2s-2(2t-2)=4t(n-1)-2(n-1)+36t^2+12t-2(s-1)$. By the induction
hypothesis, $\pi_{n}'$ is potentially $A'(3t+1)$-graphic, and hence
so is $\pi$.\ \ $\Box$

{\bf Proof of Theorem 4.2.}\ \ Let $t\ge3, n\ge20t^2+23t+5$ and
$\pi=(d_1,\ldots,d_n)\in GS_n$ with $\sigma(\pi)>4tn-2n-4t^2+2t+2$.
We only need to prove that $\pi$ is potentially $A'(3t+1)$-graphic.
We use induction on $n$. If $n=20t^2+23t+5$, by Lemma 4.2.6
$(s=20t^2+5t-1)$, then $\pi$ is potentially $A'(3t+1)$-graphic.
Assume that $n\ge20t^2+23t+6$. If $d_n\ge 2t-1$, by Lemma 4.2.5,
then $\pi$ is potentially $A'(3t+1)$-graphic. If $d_n\le 2t-2$, then
$\pi_{n}'$ satisfies
$\sigma(\pi_{n}')=\sigma(\pi)-2d_n>4tn-2n-4t^2+2t+2-2(2t-2)>4t(n-1)-2(n-1)-4t^2+2t+2$.
By the induction hypothesis, $\pi_{n}'$ is potentially
$A'(3t+1)$-graphic, and hence so is $\pi$.\ \ $\Box$

\section*{4.3. Proof of Theorem 4.3}

\hskip\parindent {\bf Lemma 4.3.1}\ \ {\it Let $t\ge 2, n\ge 18t+12$ and
$\pi=(d_1,\ldots,d_n)\in GS_n$ with $\sigma(\pi)>4tn-4t^2-2t$.
Then

(1)\ \ $d_i\ge 3t+2-\lceil\frac{i}{2}\rceil$ for $i=1,\ldots,2t-1$
and $d_{2t+1}\ge 2t+1$;

(2)\ \ $d_i\ge2(3t+3-i)$ for $i=2t+3,\ldots,3t+2.$}

{\bf Proof.}\ \ (1) If there is an even $s$ with $2\le s\le 2t-2$
such that $d_s\le 3t+2-\lceil\frac{s}{2}\rceil-1=3t+1-\frac{s}{2}$,
then
$$\begin{array}{lll}
\sigma(\pi)&\leq&(s-1)(n-1)+(3t+1-\frac{s}{2})(n-s+1)\\
&=&\frac{s^2}{2}-s(3t-\frac{n}{2}+\frac{5}{2})+3tn+3t+2.
\end{array}$$
Denote $f(s)=\frac{s^2}{2}-s(3t-\frac{n}{2}+\frac{5}{2})+3tn+3t+2$.
Since $2\le s\le 2t-2$, we have that
$$\begin{array}{lll}
\sigma(\pi)&\leq& f(s)\leq \max\{f(2),f(2t-2)\}\\
&=&\max\{3tn+n-3t-1,4tn-4t^2-n+9\}\\
&=&\max\{4tn-4t^2-2t-[(t-1)n+4t^2+t+1],4tn-4t^2-2t-(n-2t-9)\}\\
&\le&4tn-4t^2-2t,
\end{array}$$ a contradiction. Hence $d_i\ge 3t+2-\lceil\frac{i}{2}\rceil$
for even $s$ with $2\le s\le 2t-2$.

If there is an odd $s$ with $1\le s\le2t+1$ such that $d_s\le
3t+2-\lceil\frac{s}{2}\rceil-1=3t+1-\frac{s+1}{2}$, then
$$\begin{array}{lll}
\sigma(\pi)&\leq&(s-1)(n-1)+(3t+1-\frac{s+1}{2})(n-s+1)\\
&=&\frac{s^2}{2}-s(3t-\frac{n}{2}+2)+3tn+3t-\frac{n}{2}+\frac{3}{2}.
\end{array}$$
Denote
$g(s)=\frac{s^2}{2}-s(3t-\frac{n}{2}+2)+3tn+3t-\frac{n}{2}+\frac{3}{2}$.
Since $1\leq s \leq 2t+1$, we have that
$$\begin{array}{lll}
\sigma(\pi)&\leq&g(s) \leq \max\{g(1),g(2t+1)\}\\
&=&\max\{3tn,4tn-4t^2-2t\}\\
&=&\max\{4tn-4t^2-2t-[(n-4t)(t-1)+(n-6t)],4tn-4t^2-2t\}\\
&=&4tn-4t^2-2t,
\end{array}$$ a contradiction. Hence $d_i\ge 3t+2-\lceil\frac{i}{2}\rceil$
for odd $s$ with $1\le s\le2t+1$, that is, $d_i\ge
3t+2-\lceil\frac{i}{2}\rceil$ for odd $s$ with $1\le s\le 2t-1$ and
$d_{2t+1}\ge 3t+2-\lceil\frac{2t+1}{2}\rceil=2t+1.$

(2) If there is an $s$ with $2t+3\le s\le 3t+2$ such that $d_s\le
6t-2s+5$, by Theorem 2.2, then
$$\begin{array}{lll}
\sigma(\pi)&=&\sum\limits_{i=1}^{n}d_i
=\sum\limits_{i=1}^{s-1}d_i+\sum\limits_{i=s}^{n}d_i
\leq((s-2)(s-1)+\sum\limits_{i=s}^{n}\min\{s-1,d_i\})+\sum\limits_{i=s}^{n}d_i\\
&=&(s-2)(s-1)+2\sum\limits_{i=s}^{n}d_i\leq(s-2)(s-1)+2(6t-2s+5)(n-s+1)\\
&=&5s^2-(12t+4n+17)s+12tn+12t+10n+12.
\end{array}$$
Denote $h(s)=5s^2-(12t+4n+17)s+12tn+12t+10n+12$. Since $2t+3\le s\le
3t+2$, we have that
$$\begin{array}{lll}
\sigma(\pi)&\leq&h(s) \leq \max\{h(2t+3),h(3t+2)\}\\
&=&\max\{4tn-2n-4t^2+2t+6,2n+9t^2-3t-2\}\\
&=&\max\{4tn-4t^2-2t-(2n-4t-6),\\
& &4tn-4t^2-2t-[(n-4t)(4t-2)+(3t^2-7t+2)]\}\\
&\le&4tn-4t^2-2t,
\end{array}$$a contradiction.  \ \ $\Box$

We now define a new graph $G(3t+2)$ for $t\ge 2$ as follows. Let
$V(K_{2t+2})=\{v_1,\ldots,v_{2t+2}\}$ and $G(3t+2)$ be the graph
obtained from the graph $K_{2t+2}-v_{2t}v_{2t+2}$ by adding new
vertices $x_1,\ldots,x_{t}$ and joining $x_i$ to $v_1,\ldots,v_{2i}$
for $1\le i\le t$.

{\bf Lemma 4.3.2}\ \ {\it If $t\ge2$ and $G$ is any 2-tree on $3t+2$
vertices, then $G(3t+2)$ contains $G$.}

{\bf Proof.}\ \  We use induction on $t$. Let $xy\in C(G)$ so that
$xy$ is attached by the ear $z$. Denote $H=G-\{x,y,z\}$. By Lemma
4.3, $H$ is a spanning subgraph of some 2-tree $T$ on $3(t-1)+2$
vertices and $T\not=T(3(t-1)+2)$. Assume that $t=2$. Then $G(8)
-\{v_1,v_2,x_1\}=K_5-E(P_4)$. Since $K_5-E(P_4)$ and $T(5)$ are the
only two 2-trees on 5 vertices, we have that $G(8) -\{v_1,v_2,x_1\}$
contains any 2-tree $T$ on 5 vertices, where $T\not=T(5)$. Hence $H$
is a spanning subgraph of $G(8)-\{v_1,v_2,x_1\}$. Putting $x,y$ and
$z$ on $v_1,v_2$ and $x_1$ respectively, we can see that $G(8)$
contains $G$. If $t\ge 3$, then $G(3t+2)
-\{v_1,v_2,x_1\}=G(3(t-1)+2)$ contains $T$, and hence contains $H$.
Putting $x,y$ and $z$ on $v_1,v_2$ and $x_1$ respectively, we can
see that $G(3t+2)$ contains $G$.\ \ $\Box$

Let $t\ge 2,\ n\ge 18t+12$ and $\pi=(d_1,\ldots,d_n)\in GS_n$ satisfy

 (i)\ \  $d_i\ge 3t+2-\lceil\frac{i}{2}\rceil$ for $i=1,\ldots,2t-1$ and
$d_{2t+1}\ge 2t+1$, and

(ii)\ \ $d_n\ge 2t$.

We now define sequence $\pi_0,\pi_1,\ldots,\pi_{3t+2}$ as follows.
Let $\pi_0=\pi$. We define the sequence
$$\pi_1=(d_{2}^{(1)},\ldots,d_{3t+2}^{(1)},d_{3t+3}^{(1)},\ldots,d_{n}^{(1)})$$
from $\pi_0$ by deleting $d_1$, decreasing the first $d_1$ remaining
nonzero terms each by one unity, and then reordering the last
$n-(3t+2)$ terms to be non-increasing.

For $2\le i\le 3t+2$ and $i\not=2t$, we define the sequence
$$
\pi_i=(d_{i+1}^{(i)},\ldots,d_{3t+2}^{(i)},d_{3t+3}^{(i)},\ldots,d_{n}^{(i)})$$
from
$$\pi_{i-1}=(d_{i}^{(i-1)},\ldots,d_{3t+2}^{(i-1)},d_{3t+3}^{(i-1)},\ldots,d_{n}^{(i-1)})$$
by deleting $d_i^{(i-1)}$, decreasing the first
$d_i^{(i-1)}$ remaining nonzero terms each by one unity, and then
reordering the last $n-(3t+2)$ terms to be non-increasing.

For $i=2t$ and $d_{2t}\ge 2t+2$, the definition of $\pi_{2t}$
is as above. For $i=2t$ and $d_{2t}=2t+1$, we define the
sequence $$
\pi_{2t}=(d_{2t+1}^{(2t)},\ldots,d_{3t+2}^{(2t)},d_{3t+3}^{(2t)},\ldots,d_{n}^{(2t)})$$
from
$$\pi_{2t-1}=(d_{2t}^{(2t-1)},\ldots,d_{3t+2}^{(2t-1)},d_{3t+3}^{(2t-1)},\ldots,d_{n}^{(2t-1)})$$
by deleting $d_{2t}^{(2t-1)}$, decreasing the first
$d_{2t}^{(2t-1)}$ remaining nonzero terms each by one unity except
for the term $d_{2t+2}^{(2t-1)}$, and then reordering the last
$n-(3t+2)$ terms to be non-increasing.

By the definition of $\pi_{3t+2}$, the following Proposition 4.3.1 is obvious.

{\bf Proposition 4.3.1}\ \ {\it Let $t\ge 2$, $n\geq 18t+12$ and
$\pi=(d_1,\ldots,d_n)\in GS_n$ satisfy (i)--(ii). If $\pi_{3t+2}$
is graphic, then $\pi$ has a realization containing $G(3t+2)$ on
those vertices with degrees $d_1,\ldots,d_{3t+2}$.}

{\bf Lemma 4.3.3}\ \ {\it Let $t\ge 2$, $n\geq 18t+12$ and
$\pi=(d_1,\ldots,\ldots,d_n)
 \in GS_n$ satisfy (i)--(ii). For each
$\pi_i=(d_{i+1}^{(i)},\ldots,d_{3t+2}^{(i)},d_{3t+3}^{(i)},\ldots,d_{n}^{(i)})$,
let $s_i=\max\{j|d_{3t+3}^{(i)}-d_{3t+2+j}^{(i)}\le 1\}$.

(1)\ \ If $\pi$ satisfies (d) or (e), where (d) $d_1\le n-2$,
$d_{3t+2}>d_{d_2+2}$ and $d_{3t+2}-d_{d_1+2}\le 1$ and (e) $d_1\le
n-2$, $d_{3t+2}=d_{d_2+2}$ and $d_{d_2+2}=d_{d_1+2}$, then
$d_{3t+2+r}^{(3t+2)}=d_{3t+2+r}$ for $r>s_{3t+2}$.

(2)\ \ If $\pi$ satisfies (f) or (g), where (f) $d_1=n-1$, $d_2\le
n-2$ and $d_{3t+2}=d_{d_2+2}$ and (g) $d_1\le n-2$,
$d_{3t+2}=d_{d_2+2}$ and $d_{d_2+2}>d_{d_1+2}$, then
$d_{3t+2+r}^{(3t+2)}=d_{3t+2+r}^{(1)}$ for $r>s_{3t+2}$.}

{\bf Proof.}\ \ The proof of Lemma 4.3.3 is similar to that of Lemma 3.5.\ \ $\Box$

{\bf Lemma 4.3.4}\ \ {\it Let $t\ge 2$, $n\ge 18t+12$ and
$\pi=(d_1,\ldots,d_n)\in GS_n$ with $d_n\ge 2t$ and
$\sigma(\pi)>4tn-4t^2-2t$. Then $\pi$ is potentially
$A''(3t+2)$-graphic.}

{\bf Proof.}\ \ We only need to prove the following two claims.

{\bf Claim 1}\ \ {\it If $\pi$ satisfies one of (d)--(g), then $\pi$
is potentially $A''(3t+2)$-graphic.}

{\it Proof of Claim 1.}\ \ If $d_{3t+2}\ge 6t+1$, by Theorem 2.5,
then $\pi$ has a realization containing $K_{3t+2}$, implying that
$\pi$ is potentially $A''(3t+2)$-graphic by Theorem 2.4. Assume that
$d_{3t+2}\le 6t$. By Lemma 4.3.1, we have $d_i\ge
3t+2-\lceil\frac{i}{2}\rceil$ for $i=1,\ldots,2t-1$ and $d_{2t+1}\ge
2t+1$. It is enough to prove that $\pi_{3t+2}$ is graphic by Lemma
4.3.2 and Proposition 4.3.1. If $\pi$ satisfies (d) or (e), by Lemma
4.3.3(1), then
$$\pi_{3t+2}=(d_{3t+3}^{(3t+2)},\ldots,d_{3t+2+s_{3t+2}}^{(3t+2)},d_{3t+2+s_{3t+2}+1},\ldots,d_n).$$
If $\pi$ satisfies (f) or (g), by Lemma 4.3.3(2), then
$$\pi_{3t+2}=(d_{3t+3}^{(3t+2)},\ldots,d_{3t+2+s_{3t+2}}^{(3t+2)},d_{3t+2+s_{3t+2}+1}^{(1)},\ldots,d_{n}^{(1)}).$$
If $s_{3t+2}<n-(3t+2)$, by $d_{3t+3}^{(3t+2)}\le d_{3t+2}\le 6t$ and
$d_n\ge d_n^{(1)}\ge d_n-1\ge 2t-1\ge 3$, then we have that
$$\begin{array}{lll}
\frac{1}{2t-1}\left\lfloor\frac{(6t+2t-1+1)^2}{4}\right\rfloor&\le &\frac{(6t+2t-1+1)^2}{4(2t-1)}\\
&=&\frac{64t^2}{8t-4}\\
&=&\frac{8t(8t-4)+4(8t-4)+16}{8t-4}\\
&\le&8t+4+2\le n-(3t+2).
\end{array}$$
By Theorem 2.3, $\pi_{3t+2}$ is graphic. If $s_{3t+2}=n-(3t+2)$,
then $d_{3t+3}^{(3t+2)}-d_n^{(3t+2)}\le 1$. Denote $d_n^{(3t+2)}=m$.
If $m=0$, by $d_{3t+3}^{(3t+2)}\le 1$ and $\sigma(\pi_{3t+2})$ being
even, then $\pi_{3t+2}$ is clearly graphic. If $m\ge 1$, then
$d_{3t+3}^{(3t+2)}\le m+1$, and hence
$$\frac{1}{m}\left\lfloor\frac{(m+1+m+1)^2}{4}\right\rfloor=\frac{(m+1)^2}{m}\leq m+3\le 6t+3\le
n-(3t+2).$$ By Theorem 2.3, $\pi_{3t+2}$ is also graphic. This
proves Claim 1.\ \ $\Box$

{\bf Claim 2}\ \ {\it If $\pi$ satisfies one of (a)--(c), then $\pi$
is potentially $A''(3t+2)$-graphic.}

{\it Proof of Claim 2.}\ \ We use induction on $t$. Let $G$ be any
2-tree on $3t+2$ vertices, and let $x,y,z\in V(G)$ so that $xy\in
C(G)$ and $xy$ is attached by the ear $z$. Denote $H=G-\{x,y,z\}$.
By Lemma 4.3, $H$ is a spanning subgraph of some 2-tree $T$ on
$3(t-1)+2$ vertices and $T\not=T(3(t-1)+2)$.

If $t=2$, by the definition of $\rho$, we can see that $n-2\ge 46$
and $\sigma(\rho)=\sigma(\pi)-2d_1-2d_2+2>8n-20-4(n-1)+2=4(n-2)-6$.
By Theorem 2.7(4) and Theorem 2.4, $\rho$ has a realization $G_1$ in
which the subgraph $F$ induced by the vertices with degrees
$\rho_1,\ldots,\rho_5$ contains $K_5-E(P_4)$. Note that $K_5-E(P_4)$
and $T(5)$ are the only two 2-trees on 5 vertices. Denote $F'$ to be
the graph obtained from $F$ by adding three new vertices $x,y,z$
such that $x,y$ are adjacent to each vertex of $F$ and $xy,xz,yz\in
E(F')$. Since $F$ contains $T$ (and hence $H$), we have that $F'$
contains $G$. By the arbitrary of $G$, $F'$ contains every 2-trees
on 8 vertices. Moreover, by Lemma 4.4 (the case of $k=8$), we have
$\rho_1=d_3-2,\rho_2=d_4-2,\ldots,\rho_6=d_{8}-2$. This implies that
$\pi$ has a realization $G'$ in which the subgraph induced by the
vertices with degrees $d_1,\ldots,d_8$ contains $F'$. In order
words, $\pi$ is potentially $A''(8)$-graphic.

If $t\ge 3$, by the definition of $\rho$, then $n-2\ge 18(t-1)+12$,
$\rho_{n-2}\ge 2t-2=2(t-1)$ and
$\sigma(\rho)=\sigma(\pi)-2d_1-2d_2+2>4tn-4t^2-2t-4(n-1)+2
=4(t-1)(n-2)-4(t-1)^2-2(t-1)$. If $\rho$ satisfies one of (d)--(g),
by Claim 1, then $\rho$ is potentially $A''(3(t-1)+2)$-graphic. If
$\rho$ satisfies one of (a)--(c), by the induction hypothesis, then
$\rho$ is also potentially $A''(3(t-1)+2)$-graphic. This implies
that $\rho$ has a realization $G_1$ in which the subgraph $F$
induced by the vertices with degrees $\rho_1,\ldots,\rho_{3(t-1)+2}$
contains every 2-tree on $3(t-1)+2$ vertices. Denote $F'$ to be the
graph obtained from $F$ by adding three new vertices $x,y,z$ such
that $x,y$ are adjacent to each vertex of $F$ and $xy,xz,yz\in
E(F')$. Since $F$ contains $T$ (and hence $H$), we have that $F'$
contains $G$. By the arbitrary of $G$, $F'$ contains every 2-trees
on $3t+2$ vertices. Moreover, by Lemma 4.4 (the case of $k=3t+2$),
we have $\rho_1=d_3-2,\rho_2=d_4-2,\ldots,\rho_{3t}=d_{3t+2}-2$.
This implies that $\pi$ has a realization $G'$ in which the subgraph
induced by the vertices with degrees $d_1,\ldots,d_{3t}$ contains
$F'$. In order words, $\pi$ is potentially $A''(3t)$-graphic. This
proves Claim 2.\ \ $\Box$

{\bf Lemma 4.3.5}\ \ {\it Let $t\ge 2$, $n=18t+12+s$ with $0\le s\le 20t^2+13t$ and
$\pi=(d_1,\ldots,d_n)\in GS_n$ with $\sigma(\pi)>4tn+36t^2+24t-2s$. Then $\pi$ is potentially
$A'(3t+2)$-graphic.}

{\bf Proof.}\ \ We use induction on $s$. If $s=0$, then $n=18t+12$
and $\sigma(\pi)\ge 4tn+36t^2+24t+2=2n(3t+2-2)+2$. By Theorem 2.6,
$\pi$ has a realization containing $K_{3t+2}$. and hence $\pi$
is potentially $A'(3t+2)$-graphic. Suppose now that $1\le s\le
 20t^2+13t$. Then $\sigma(\pi)>4tn-4t^2-2t$. If $d_n\ge 2t$,
then $\pi$ is potentially $A'(3t+2)$-graphic by Lemma 4.3.4. If
$d_n\le 2t-1$, then $\pi_{n}'=(d_{1}',\ldots,d_{n-1}')$ satisfies
$\sigma(\pi_{n}')=\sigma(\pi)-2d_n>4tn+36t^2+24t-2s-2(2t-1)=4t(n-1)+36t^2+24t-2(s-1)$.
By the induction hypothesis, $\pi_{n}'$ is potentially
$A'(3t+2)$-graphic, and hence so is $\pi$.\ \ $\Box$

{\bf Proof of Theorem 4.3.}\ \  Let $t\ge 2$, $n\ge 20t^2+31t+12$
and $\pi=(d_1,\ldots,d_n)\in GS_n$ with $\sigma(\pi)>4tn-4t^2-2t$.
We only need to prove that $\pi$ is potentially $A'(3t+2)$-graphic.
We use induction on $n$. If $n=20t^2+31t+12$, by Lemma 4.3.5
($s=20t^2+13t$), then $\pi$ is potentially $A'(3t+2)$-graphic.
Assume that $n\ge 20t^2+31t+13$. If $d_n\ge 2t$, by Lemma 4.3.4,
then $\pi$ is potentially $A'(3t+2)$-graphic. If $d_n\le 2t-1$, then
$\pi_{n}'$ satisfies $\sigma(\pi_{n}')=\sigma(\pi)-2d_n>
4tn-4t^2-2t-2(2t-1)>4t(n-1)-4t^2-2t$. By the induction hypothesis,
$\pi_{n}'$ is potentially $A'(3t+2)$-graphic, and hence so is
$\pi$.\ \ $\Box$

\vskip 0.2cm

\noindent{\bf References}

\bref{[1]}{J.A. Bondy and U.S.R. Murty,}{ Graph Theory with
Applications,}{} { The Macmillan Press, London, 1976.}

\vskip 0.1cm

\bref{[2]}{P. Bose, V. Dujmovic, D. Krizanc, S. Langerman, P. Morin,
D.R. Wood and S. Wuhrer, A characterization of the degree sequences
of 2-trees,}{ J. Graph Theory,}{ 58}{ (2008), 191--209.}

\vskip 0.1cm

\bref{[3]}{L.Z. Cai, On spanning 2-trees in a graph,}{ Discrete
Appl. Math.,}{ 74}{ (1997), 203-216.}

\vskip 0.1cm

\bref{[4]} {P. Erd\H{o}s and T. Gallai, Graphs with prescribed
degrees of vertices (Hungarian),}{ Mat. Lapok,} { 11}{ (1960),
264--274.}

\vskip 0.1cm

\bref{[5]}{P. Erd\H{o}s, M.S. Jacobson and J. Lehel, Graphs
realizing the same degree sequences and their respective clique
numbers, in: Y. Alavi et al., (Eds.),}{ Graph Theory, Combinatorics
and Applications,}{}{ Vol.1, John Wiley \& Sons, New York, 1991,
439--449.}

\vskip 0.1cm

\bref{[6]}{R.J. Gould, M.S. Jacobson and J. Lehel, Potentially
$G$-graphical degree sequences, in: Y. Alavi et al., (Eds.),}{
Combinatorics, Graph Theory, and Algorithms,}{}{ Vol.1, New Issues
Press, Kalamazoo Michigan, 1999, 451--460.}

\vskip 0.1cm

\bref{[7]}{D.J. Kleitman and D.L. Wang, Algorithm for constructing
graphs and digraphs with given valences and factors,}{ Discrete
Math.,}{ 6}{ (1973), 79--88.}

\vskip 0.1cm

\bref{[8]}{C.H. Lai, A note on potentially $K_4-e$ graphical
sequences,}{ Australas. J. Combin.,}{ 24}{ (2001), 123--127.}

\vskip 0.1cm

\bref{[9]}{C.H. Lai and L.L. Hu,  An extremal problem on potentially
$K_{r+1}- H$-graphic sequences}{ Ars Combin.,}{ 94}{ (2010),
289--298.}

\vskip 0.1cm

\bref{[10]}{C.H. Lai and L.L. Hu, Potentially $K_m-G$-graphical
sequences: a survey,}{ Czechoslovak Math. J.,}{ 59(134)}{ (2009),
1059--1075.}

\vskip 0.1cm

\bref{[11]}{J.S. Li and Z.X. Song, On the potentially $P_k$-graphic
sequence,} { Discrete Math.,}{ 195} { (1999), 255--262.}

\vskip 0.1cm

\bref{[12]}{J.H. Yin and J.S. Li, Two sufficient conditions for a
graphic sequence to have a realization with prescribed clique
size,}{ Discrete Math.,}{ 301}{ (2005), 218--227.}

\vskip 0.1cm

\bref{[13]}{J.H. Yin and J.S. Li, An extremal problem on potentially
$K_{r,s}$-graphic sequences,}{ Discrete Math.,}{ 260}{ (2003),
295--305.}

\vskip 0.1cm

\bref{[14]}{J.H. Yin and J.S. Li, A variation of a conjecture due to
 Erd\H{o}s and S\'{o}s,}{ Acta Math. Sin. Engl. Ser.,}{ 25}{ (2009), 795--802.}

 \vskip 0.1cm

\bref{[15]}{D.Y. Zeng and J.H. Yin, An extremal problem for a
graphic sequence to have a realization containing every $2$-tree
with prescribed size,}{ Discrete Math. Theor. Comput. Sci.,}{
17(3)}{ (2016), 315--326.}

\end{document}